\documentclass[12pt]{article}
\usepackage[dvipdfmx]{graphicx,color}
\usepackage{amsmath,amsfonts,amssymb}
\usepackage{hyperref}
\usepackage{ulem}

\usepackage[numbers]{natbib}
\newcommand{\ds}{\displaystyle }

\newcommand{\vc}[1]{{\boldsymbol #1}}

\newcommand{\sr}[1]{{\cal #1}}
\newcommand{\dd}[1]{\mathbb{#1}}

\newcommand{\eq}[1]{(\ref{eq:#1})}
\newcommand{\lem}[1]{Lemma~\ref{lem:#1}}
\newcommand{\cor}[1]{Corollary~\ref{cor:#1}}
\newcommand{\thr}[1]{Theorem~\ref{thr:#1}}

\newcommand{\rem}[1]{Remark~\ref{rem:#1}}
\newcommand{\fig}[1]{Figure~\ref{fig:#1}}
\newcommand{\app}[1]{Appendix~\ref{app:#1}}
\newcommand{\sectn}[1]{Section~\ref{sect:#1}}

\newcommand{\sect}[1]{\ref{sect:#1}}

\newcommand{\ol}{\overline}

\newcommand{\pend}{\hfill \thicklines \framebox(6.6,6.6)[l]{}}

\newenvironment{proof}{\noindent {\sc  Proof.} \rm}{\pend}
\newenvironment{proof*}[1]{\noindent {\sc  #1} \rm}{\pend}

\newtheorem{theorem}{Theorem}[section]
\newtheorem{lemma}{Lemma}[section]

\newtheorem{remark}{Remark}[section]

\newtheorem{corollary}{Corollary}[section]

\newcommand{\setsection}[2] {
\setcounter{section}{#1}
\setcounter{subsection}{0}
\setcounter{equation}{0}
\setcounter{conjecture}{0}
\setcounter{assumption}{0}
\setcounter{question}{0}
\setcounter{definition}{0}
\setcounter{theorem}{0}
\setcounter{corollary}{0}
\setcounter{lemma}{0}
\setcounter{proposition}{0}
\setcounter{remark}{0}
\setcounter{appen}{0}
\setsection*{\large \bf \thesection. #2}}

\newcommand{\setnewcounter} {
\setcounter{subsection}{0}
\setcounter{equation}{0}
\setcounter{conjecture}{0}
\setcounter{assumption}{0}
\setcounter{question}{0}
\setcounter{definition}{0}
\setcounter{theorem}{0}
\setcounter{corollary}{0}
\setcounter{lemma}{0}
\setcounter{proposition}{0}
\setcounter{remark}{0}
}

\begin{document}
\title{\bf \Large Customer sojourn time in $GI/GI/1$ feedback queue\\ in the presence of heavy tails}

\author{Sergey Foss\\Heriot-Watt University and \\ Novosibirsk State University\footnote{Research of S.~Foss is supported
by RSF research grant No. 17-11-01173. He also thanks EPFL, Lausanne for their hospitality}\\ \and Masakiyo Miyazawa\\ Tokyo University of Science\footnote{Research of M.~Miyazawa is supported in part by
JSPS KAKENHI Grant No. JP16H02786}}
\date{June 17, 2018 (to appear in the Journal of Statistical Physics)}

\maketitle

\begin{abstract}
We consider a single-server $GI/GI/1$ queueing system with feedback. 
We assume the service time distribution to be (intermediate) regularly varying. 
We find the tail asymptotics for a customer's sojourn time in two cases: the customer arrives in an empty system, and the customer arrives in the system
in the stationary regime. In particular, in the case of Poisson input we 
obtain more explicit formulae than those in the general case. As auxiliary results, we find the tail asymptotics for the busy period
distribution in a single-server queue with an intermediate varying service times distribution and establish the principle-of-a-single-big-jump equivalences that characterise the asymptotics.

\end{abstract}
  
\begin{quotation}
\end{quotation}

\section{Introduction}
\label{sect:introduction}

In queueing theory, the sojourn time $U$ of a customer in a queueing system is one of important characteristics, this is the time from its arrival instant to departure instant. 
In general, the distribution of $U$ is hard to find analytically, and research interest is directed to the asymptotics of the tail probability, 
$\dd{P}(U > x)$, as $x \to \infty$ under various stochastic assumptions. Among them, the following assumption is typically used.
\begin{itemize}

\item [(\sect{introduction})] Delaying arrivals leads to increase of the characteristic. 

\end{itemize}
The assumption (\sect{introduction}) is known as a {\it monotonicity} property, it plays an important role in the asymptotic analysis of various characteristics. Another important factor for the tail asymptotic problem is the heaviness of service time distributions. A nonnegative random variable $X$ is said to have a heavy tail distribution if $\dd{E} \exp (sX) = \infty$ for all $s > 0$, and a light tail distribution, otherwise. If service time distributions have heavy tails (or light tails), we talk about a heavy (or light) tail regime of the system.

Under the heavy tail regime, the sojourn time problem is relatively well studied for the system satisfying the monotonicity assumption  (\sect{introduction}) (see e.g. \cite{BaFo2004}, \cite{JeMoZw2004}, \cite{FoKo2012}, \cite{FoMy2014} and references therein).
However, if a customer separately takes multiple services, then the monotonicity is violated, in general. This is a common phenomenon in queueing networks. Clearly, {\it non-monotone} characteristics are also important, but we are unaware of any asymptotic results for them in the presence of heavy tails.
So we challenge the tail asymptotics problem of a non-monotone characteristic under the heavy tail regime, using a relatively simple model as one of the first attempts.

We consider the following single server system. Exogenous customers arrive at the system subject to $i.i.d.$ inter-arrival times $\{t_n\}$ with finite mean $a$, join at the end of a queue, and are served in the first-in first-out order by a single server with $i.i.d.$ service times $\{\sigma_n\}$. When a customer completes service, it returns to the end of the queue for another service, with probability $p\in (0,1)$, or leaves the system, with probability $q=1-p$. Both transitions are independent of everything else. Note that if a customer requires more than one service, its sojourn time is the sum of several periods of waiting in queue with following service. One can see that the sojourn time $U$ does not have natural monotonicity properties with respect to inter-arrival times. We refer this system as a $GI/GI/1$ feedback queue. If $p=0$ (no customer returns), then it is reduced to a standard $GI/GI/1$ queue (e.g., see \cite{As2003}).

Because of independent feedback of customers completing service, the $i$-th arriving customer, customer $i$ for short, requires a geometric number of services, say $K_i$, and their total duration is $\sum_{j=1}^{K_i} \sigma_i^{(j)}$ where $\{\sigma_i^{(j)}\}$ are $i.i.d.$ with finite mean $b$ and do not depend on $K_{i}$. Throughout the paper, we assume that the system is stable, that is, $\rho \equiv \lambda b/q < 1$, where $\lambda = 1/a$, and the service-time distribution is heavy-tailed. More precisely, we
assume the distribution to be {\it intermediate regularly varying} -- see the list of heavy-tailed distributions at the end of this section and in  \app{properties}. 

For this $GI/GI/1$ feedback queue, we derive the asymptotics of the  probability $\dd{P}(U>x)$ for the sojourn time $U$ of customer $1$ as $x\to\infty$ under the heavy tail regime. This $U$ depends on the state of the system just before the arrival instant of customer $1$. We consider two scenarios for the system state found by customer $1$. The first one is that customer $1$ enters an empty queue. It is the most difficult because
coefficients appearing in the asymptotics are sensitive not only to the first moments, but to the whole distribution of inter-arrival and service times, in general. The second one is that customer $1$ enters a stationary queue. In both cases, we also obtain simpler formulae in the case when the input is Poisson. In particular, the first moment of the sojourn time can be obtained in a closed form in this special case. A part of this information will be used for the general renewal arrival case. 

Our model may be considered as a particular case of a two-server {\it generalised Jackson network} where server $i=1,2$ has service times  $\{\sigma_n^{(i)}\}$. Each family of service times is $i.i.d.$ and they are mutually independent. 
Customers arrive in a renewal input to server 1 and join the queue there. After service completion at server 1, a customer either leaves the network,  with probability $p_{10}$, or joins the queue to the second server, with probability $p_{12}=1-p_{10}$. Similarly, after service completion at server 2, a customer either leaves the
network, with probability $p_{20}$, or joins the queue to server 1, with probability $p_{21}=1-p_{20}$. Customers are server in the order of their (external and internal) arrival
to the servers. 

If we let $\sigma_n^{(2)}\equiv 0$ and $p_{12}p_{21}=p$, we obtain our
model as a particular case indeed. So the study of our model is not only of interest itself, but also opens a window to analysing a broad class of more general models.

In the $GI/GI/1$ feedback queue, one may change the service order in such a way that each customer continuously gets service without interruption when it completes service and returns to the queue. Then such a system is nothing else than the standard $GI/GI/1$ queue with ``new'' $i.i.d.$ service times $\sum_{j=1}^{K_i} \sigma_i^{(j)}$, and the sojourn time is again the sum of the waiting time and of the (new) service time. Although the busy period, which is the time from the moment when the system becomes non-empty to the moment when it is again empty, is unchanged by this modification, the sojourn time does change. We will use this modified system to study the tail asymptotic of the busy period.

Thus, our analysis is connected to the standard $GI/GI/1$ queue. In this case, there is no feedback, and the monotonicity is satisfied. Hence, the waiting time of a tagged customer is a key characteristic because the sojourn time $U$ is the sum of the waiting and service times, which are independent. In particular, the stationary waiting time is a major target for the tail asymptotic analysis. Let $u_{0}$ be the unfinished work found by the ``initial'' customer $1$ that arrives at the system at time $0$, and let $W_{n}$ be the waiting times of the $n$-th arriving customer. Then, $W_{1} = u_{0}$ and we have the {\it Lindley recursion}:
\begin{align}\label{eq:lindley}
W_{n+1}= \max (0, W_n+\sigma_{n}-t_n), \qquad n \ge 1.
\end{align}
We assume both inter-arrival and service times to have finite means, $a=\dd{E}t_n$ and $b=\dd{E}\sigma_n$. 
Here $W_n$ forms a Markov chain which is {\it stable} (i.e. converges in distribution to the limiting/stationary random variable $W=W_{\infty}$) if the {\it traffic intensity} $\rho :=b/a $ is less than 1. It is well-known (see e.g. \cite{As2003}) that if $u_0=0$, then $W$ coincides in distribution
with the supremum $M=\sup_{n\ge 0} \sum_{i=1}^n (\sigma_n-t_n)$ of a random walk with increments $\sigma_n-t_n$.

The tail asymptotics for $\dd{P}(M>x)$ as $x\to\infty$ is known under the {\it light-tail} and {\it heavy-tail} regimes. In the case of light tails, there are three types of the tail asymptotics, depending on properties
of the moment generating function $\varphi (s) = \dd{E} \exp (s\sigma)$ -- see e.g. \cite{FoPu(2011)} and references therein. In the case of heavy tails, the tail asymptotics
are known in the class of so-called {\it subexponential distributions} and are  based on the {\it principle of a single big jump (PSBJ)}: 
$M$ takes a large value if one of the service times is large. This PSBJ has been used for the asymptotic analysis in several other (relatively simple) stable queueing models,  for
a number of characteristics (waiting time, sojourn time, queue length, busy cycle/period, maximal data, etc.) that possess the monotonicity property  (\sect{introduction}) (see e.g. \cite{BaFo2004}). 

Our proofs rely on the tail asymptotics for the first and stationary busy periods of the system. We establish the PSBJ for the busy period
first. This allows us to establish the principle for the sojourn time since the tail distribution asymptotics of the busy
period is of the same order with that of the sojourn time. Then insensitivity properties of the {\it intermediate varying  distributions} (see Appendix A again) allow us to compute the exact tail asymptotics for the sojourn time. The main result from \cite{FoZa2003} is a key tool
in our analysis.  

The paper is organised as follows. \sectn{main} formally introduces the model and presents main results. \sectn{psbj} states the tail asymptotic of the busy period and the PSBJ. All theorems from Sections \sect{main} and \sect{psbj} are proved in \sectn{proofs}. The Appendix consists of three parts. Part A contains an overview on basic properties of heavy-tailed distributions and part B the proof of Corollary 3.1. In part C, we propose an alternative
approach to the proof of Corollary 3.2.

Throughout the paper, we use the following notation: $1(\cdot)$ is the indicator function of the event ``$\cdot$''. For two positive functions $f$ and $g$, we write $f(x)\sim g(x)$ if $f(x)/g(x)\to 1$ as $x\to\infty$, $f(x) \gtrsim g(x)$ if
$\liminf_{x\to\infty} f(x)/g(x) \ge 1$ and $f(x)\lesssim g(x)$ if
$\limsup_{x\to\infty} f(x)/g(x) \le 1$. For a distribution function $F$, its tail $\ol{F}$ is defined as
\begin{align*}
  \ol{F}(x) = 1 - F(x).
\end{align*}
For random variables $X, Y$ with distributions $F, G$, respectively, $X =_{st} Y$ if $F = G$, and $X \le_{st} Y$ if $\ol{F}(x) \le \ol{G}(x)$ for all real $x$. Two families of events $A_x$ and $B_x$ of non-zero probabilities are equivalent, $A_x \simeq B_x$, if
\begin{align}\label{strongeq}
  \dd{P} (A_x\Delta B_x) = o (\dd{P}(A_x)), \quad \mbox{as} \quad x\to\infty,
\end{align}
where $A_x\Delta B_x= (A_x\setminus B_x)\cup (B_x\setminus A_x)$ is the symmetric difference of $A_x$ and $B_x$. Note that equivalence $A_x\simeq B_x$ is symmetric, since
\begin{align*}
 |\dd{P}(B_x) - \dd{P}(A_x)| \le \dd{P}(A_x \Delta B_x).
\end{align*}
Note also that $A_x\simeq B_x$ is stronger than equivalence $\dd{P}(A_{x}) \sim \dd{P}(B_{x})$.

We complete the Introduction by a short\\ 
\noindent {\large \it Summary of main classes of heavy tail distributions} \smallskip\\
\indent In this paper, we are concerned with several classes of heavy tail distributions. We list their definitions below. Their basic properties are discussed in \app{properties}.

In all definitions below, we assume that $\overline{F}(x)>0$ for all $x$.

\begin{enumerate}
\item Distribution $F$ on the real line belongs to the class ${\cal L}$ of {\it long-tailed} distributions if, 
for some $y>0$ and 
as $x\to\infty$,
\begin{equation}
\label{eq:long}
\frac{\overline{F}(x+y)}{\overline{F}(x)} \rightarrow 1
\end{equation}
(we may write equivalently $\overline{F}(x+y) \sim \overline{F}(x)$).

\item Distribution $F$ on the positive half-line belongs to the class ${\cal S}$ of {\it subexponential} distributions if
\begin{align*}
  \int_0^x F(dt) \overline{F}(x-t) \sim \overline{F}(x) 
\quad \mbox{as} \quad x\to\infty ,
\end{align*}
which is equivalent to $\int_{x}^{\infty} F(dt) F(x-t) \sim 2 \ol{F}(x)$. Distribution $F$ of a real-valued random variable $\xi$ is subexponential if distribution $F^+(x) = F(x) 1(x\ge 0)$ of random variable $\xi^+ =\max (0,\xi )$ is subexponential.

\item Distribution $F$ on the real line belongs to the class ${\cal S}^{*}$ of {\it strong subexponential} distributions if
$m^+(F) \equiv \int_{0}^{\infty} \ol{F}(x) dx$ is finite and
$$
\int_0^x \overline{F}(y) \overline{F}(x-y) dy \sim
2m^+(F) \overline{F}(x) \quad \mbox{as} \quad x\to\infty .
$$

\item Distribution $F$ on the real line belongs to the class ${\cal D}$ of {\it dominantly varying} distributions  if there exists
$\alpha >1$ (or, equivalently, for all $\alpha >1$) such that
\begin{align}
\label{eq:DV}
  \liminf_{x\to\infty}
\frac{\overline{F}(\alpha x)}{\overline{F}(x)} > 0
\end{align}

\item Distribution $F$ on the real line belongs to the class ${\cal IRV}$ of {\it intermediate regularly varying} distributions if
\begin{align}
\label{eq:IRV}
  \lim_{\alpha \downarrow 1} \liminf_{x\to\infty}
\frac{\overline{F}(\alpha x)}{
\overline{F}(x)} = 1.
\end{align}

\item Distribution $F$ on the real line belongs to the class ${\cal RV}$ of {\it regularly varying} distributions if, for some 
$\beta >0$, 
\begin{align}
\label{eq:RV}
  \overline{F}(x) = x^{-\beta} L(x),
\end{align}
where $L(x)$ is a {\it slowly varying} function, i.e. $L(cx)\sim L(x)$ as
$x\to\infty$, for any $c>0$.

\end{enumerate}

The following relations between the classes introduced above may be found,
say, in the books \cite{EKM1997} or 
\cite{FKZ2013}: in the class of distributions $F$ with finite $m^+(F)$, 
\begin{align}
\label{eq:class order}
  {\cal RV} \subset {\cal IRV} \subset {\cal L\cap D} \subset {\cal S^*} \subset {\cal S} \subset {\cal L}.
\end{align}

\section{The modelling assumptions and main results}
\label{sect:main}
\setnewcounter

In this section, we describe the dynamics of the sojourn time of a tagged customer (customer $1$), and present main results on the tail asymptotics of its sojourn time. In Subsection 2.1, we formally introduce the $GI/GI/1$ feedback queue and then, in Subsection 2.2, a particular $M/GI/1$ feedback queue, which has Poisson arrivals. Then the main results are presented in Subsection 2.3.

\subsection{$GI/GI/1$ feedback queue}
\label{sect:GI/GI/1}

Let $K$ be the number of services of the tagged customer until its departure. By the feedback assumption, $K$ is geometrically distributed with parameter $p$, that is,
\begin{align}
\label{eq:K 1}
  \dd{P}(K=k) = q p^{k-1}, \qquad k=1,2,\ldots,
\end{align}
and independent of everything else. Throughout the paper, we make the following assumptions: 
\begin{itemize}
\item [(i)] The exogenous arrival process is a renewal process with a finite mean interarrival time $a > 0$.
\item [(ii)] All the service times that start after time $0$ are $i.i.d.$ with finite mean $b > 0$, they are jointly independent of the arrival process. 
\item [(iii)] The system is stable, that is, 
\begin{align}
\label{eq:stability 1}
   \rho \equiv \lambda b/q < 1,
\end{align}
where $\lambda = 1/a$.
\end{itemize}

We denote the counting process of the exogenous arrivals by $N^{e}(\cdot) \equiv \{N^{e}(t); t \ge 0\}$. We use the notation $G$ for the service time distribution, and use $\sigma$ for a random variable subject to $G$. 

Let $(X_{0},R^{s}_{0})$ be the pair of the number of earlier customers and the remaining service time of a customer being served at time 0, where $R^{s}_{0} = 0$ if there is no customer in the system. Let $u_0$
be the waiting time of the tagged customer before the start of its first
service, Then 
\begin{align}\label{u0}
  u_{0} = R^{s}_{0} + \sum_{i=1}^{X_{0}-1} \sigma_{0,i+1},
\end{align}
where $\sigma_{0,i}$'s for $i \ge 2$ are $i.i.d.$ random variables each of which has the same distribution as $\sigma$. There are two typical scenarios for the initial distribution, that is, the distribution of $(X_{0},R^{s}_{0})$.
\begin{itemize}
\item [(\sect{main}a)] A tagged arriving customer finds the system empty. That is, $(X_{0},R^{s}_{0}) = (0,0)$.
\item [(\sect{main}b)] A tagged arriving customer finds $X_{0}$ customers and the remaining service time $R^{s}_{0}$ of the customer being served. Thus, the initial state $(X_{0},R^{s}_{0}) \ne (0,0)$. 
\end{itemize}
In this paper, we assume that the service time distribution is heavy tailed, and mainly consider the tail asymptotic of the sojourn time distribution of the $GI/GI/1$ feedback queue under the scenario (\sect{main}a). The case (\sect{main}b) when $X_{0}$ and $R_0^s$ are bounded by a constant may be studied very similarly to the case (\sect{main}a), therefore we do not analyse it. We consider the case (\sect{main}b) when $(X_{0},R^{s}_{0})$ is subject to the stationary distribution embedded at the arrival instants.

For given $(X_{0},R^{s}_{0})$, we have defined $u_{0}$. Let $X_{k}$ be the queue length behind the tagged customer when it finished its $k$th service for $k \ge 1$ when the tagged customer gets service at least $k$ times. Similarly, let $U_{k}$ be the sojourn time of the tagged customer measured from its $(k-1)$th service completion to its $k$th service completion, and let $T_{k}$ be the sojourn time of the tagged customer just after its $k$th service completion. 

We now formally define random variables $X_{k}$, $U_{k}$ and $T_{k}$ by induction. 
Let $T_0=0$. Denote the $k$th service time of the tagged customer by $\sigma_{k,0}$, while $\sigma_{k,i}, i=1,\ldots, X_{k-1}$ are the service times of the customers waiting before the tagged one on its $k$th return. Note that $\sigma_{k,i}$'s for $k \ge 1, i \ge 0$ are $i.i.d.$ random variables subject to the same distribution as $\sigma$. Then, $X_{k}$, $U_{k}$ and $T_{k}$ for $k \ge 1$ are defined as
\begin{align}
\label{eq:U k}
 & U_{k} = \begin{cases}
\sigma_{1,0} +u_{0}& k = 1, \\
 \sum_{i=0}^{X_{k-1}} \sigma_{k,i} & k \ge 2,
\end{cases},\\
 \label{eq:t k}
 & T_k = T_{k-1} + U_{k},\\
\label{eq:X k}
 & X_{k} = N^{e}(T_k) - N^{e}(T_{k-1}) + N^{B}_{k}(X_{k-1}),
\end{align}
where $u_0$ is given by \eqref{u0}, and $N^{B}_{k}(n)$'s are $i.i.d.$ random variables each of which is subject to the Binomial distribution with parameters $n, p$. The dynamics of the sojourn time is depicted below when $X_{0} = 0$, that is, a tagged customer finds the system empty.
\begin{figure}[h] 
   \centering
   \includegraphics[height=5.5cm]{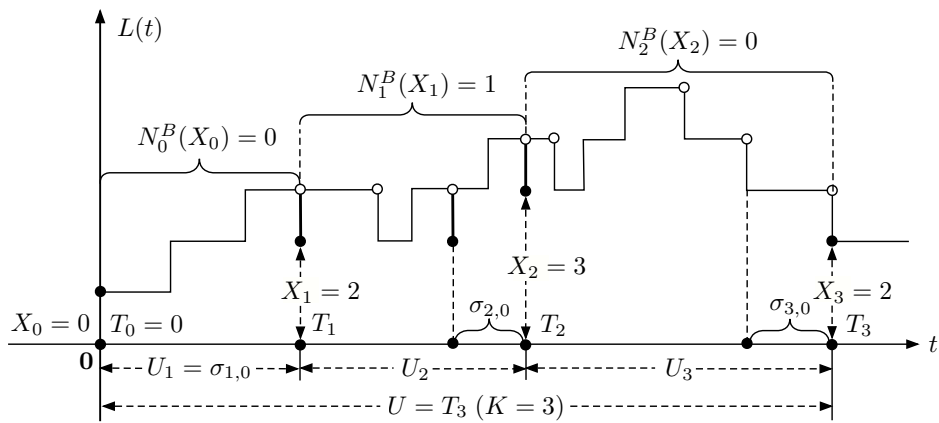}
\caption{Sample path of the queue length process $L(t)$ and $(X_{k}, U_{k})$}
\label{fig:Feedback_dynamics_1}
\end{figure}

To make clear the dependence of $X_{k}, U_{k}, T_{k}$, we introduce a filtration $\{\sr{F}_{t}; t \ge 0\}$ as
\begin{align*}
  \sr{F}_{t} = \mbox{$\sigma$-field generated by } \{X_{0}, R^{s}_{0}, (N^{e}(u), N^{s}(u), N^{r}(u)), u \le t\},
\end{align*}
where $N^{s}(t)$ and $N^{r}(t)$ are the numbers of customers who completed service and who return to the queue, respectively, up to time $t$. Clearly, $T_{k}$ is a $\sr{F}_{t}$-stopping time, and $X_{k}$ and $U_{k}$ are $\sr{F}_{T_{k}}$-measurable. Furthermore, $\sigma_{k,0}$ and $\sigma_{k,i}$ for $i \ge 1$ are independent of $\sr{F}_{T_{k-1}}$. Then $U$, the sojourn time of the tagged customer,
 may be represented as 
\begin{align}
\label{eq:U 1}
  U = u_{0} + \sum_{k=1}^{K} U_{k} = u_{0} + \sum_{k=1}^{K} \sum_{i=0}^{X_{k-1}} \sigma_{k, i}.
\end{align}

For $k\ge 0$, let $Y_{k} = \sum_{\ell = 0}^{k} X_{\ell}$ for $k \ge 0$, which is the total number of external and internal arrivals to the queue up to time $T_k$ plus the number of customers in system at time $0$. Then 
\begin{align}
\label{eq:Y recursive 1}
  Y_{k} & = X_{0} + N^{e}(T_{k}) + \sum_{\ell=1}^{k} N^{B}_{\ell}\Big(X_{\ell-1}\Big) =_{st} X_{0} + N^{e}(T_{k}) + N^{B}_{k}\Big(Y_{k-1}\Big).
\end{align}
Hence, under scenario (\sect{main}a), we have $u_{0}= X_{0} = 0$, so
\begin{align}
\label{eq:U 1 (a)}
  U = \sum_{k=1}^{K} \sum_{i=0}^{X_{k-1}} \sigma_{k i} =_{st} \sum_{i=1}^{K+Y_{K-1}} \sigma_{i},
\end{align}
while, under scenario (\sect{main}b), 
\begin{align}
\label{eq:U 1 (b)}
  U =_{st} u_{0}+ \sum_{i=1}^{K+Y_{K-1}} \sigma_{i},
\end{align}
where $\sigma_{i}$'s are $i.i.d.$ random variables each of which has the same distribution as $\sigma$. Note that $K+Y_{K-1}$ is ${\cal F}_{T_{K-1}}$-measurable that depends, in general, on all
$\sigma_i$'s of customers who arrive before $T_{K-1}$. 
This causes considerable difficulty in the asymptotic analysis of $U$.

Thus, we need to consider dependence structure in the representation of $U$. Furthermore, $\{(U_{k}, X_{k}); k \ge 0\}$ is generally not a Markov chain for a general renewal process.

On the other hand, if the arrival process $N^{e}(\cdot)$ is Poisson, then not only $\{(U_{k},X_{k}); k \ge 0\}$ but also $\{X_{k}; k \ge 0\}$ is a Markov chain with respect to the filtration $\{\sr{F}_{T_{k}}; k \ge 0\}$. In this case, we may obtain exact expressions for $\dd{E}X_k$ and then an explicit form for the tail asymptotics. 

\subsection{$M/GI/1$ feedback queue and branching process}
\label{sect:M/GI/1}

In this subsection, we assume that the exogenous arrival process is Poisson with rate $\lambda > 0$. This model is analytically studied using Laplace transforms in \cite{Taka1963}, but no asymptotic results are given there.  Note that we may consider $\{X_{k}; k\ge 0\}$ as a branching process and directly compute $\dd{E}(X_{k})$, which then will be used for the general renewal input case.

Since the Poisson process $N^{e}(\cdot)$ has independent increments, \eq{X k} is simplified to 
\begin{align}
\label{eq:X recursive 1}
 X_{k} & = N^{e}_{k}(U_{k}) + N^{B}_{k}(X_{k-1}), \qquad 1 \le k \le K,
\end{align}
using independent Poisson processes $N^{e}_{k}$ and independent Binomial  random variables $N^{B}_{k}(n)$. Furthermore, \eq{X recursive 1} can be written as
\begin{align}
\label{eq:X recursive 2}
 & X_{k} = \begin{cases}
 N_{1,1}^{e}(\sigma_{1,0} +u_{0}) + N^{B}_{1}(X_{0}), &  k=1, \\
 N_{k,0}^{e}(\sigma_{k,0}) + \sum_{i=1}^{X_{k-1}} (N_{k,i}^{e}(\sigma_{k,i}) + N^{B}_{k,i}(1)), & k \ge 2,
\end{cases}
\end{align}
where $N_{k,i}^{e}(\cdot)$'s are independent Poisson processes with rate $\lambda$. Hence, $\{X_{k}; k \ge 1\}$ is a branching process with immigration.

Due to the branching structure, we can compute the moments of $X_{k}$ explicitly. We are particularly interested in their means. 
From \eq{X recursive 2}, we have
\begin{align*}
 & \dd{E}(X_{k}) = \begin{cases}
 \lambda (b + \dd{E}(u_{0})) + p \dd{E}(X_{0}), \quad & k = 1,\\
 \lambda b + \dd{E}(X_{k-1}) r, \qquad & k \ge 2, 
\end{cases}
\end{align*}
where $r = \lambda b + p$. By the stability condition \eq{stability 1}, $r < 1$, and we have
\begin{align}
\label{eq:X k mean}
  \dd{E}(X_{k}) = \frac {1 - r^{k-1}} {1 - r} \lambda b + \dd{E}(X_{1}) r^{k-1}, \qquad k \ge 1.
\end{align}
Hence, we have a uniform bound:
\begin{align}
\label{eq:Xk bound 2}
  \dd{E}(X_{k}) \le \frac {1} {1 - r} \lambda b + \dd{E}(X_{1}), \qquad k \ge 1.
\end{align}
Furthermore, we have
\begin{align}
\label{eq:Z mean 1}
  \dd{E}(K + Y_{K-1}) &= \dd{E}(K) + \sum_{k=1}^{\infty} \sum_{\ell=0}^{k-1} \dd{E}(X_{\ell}) \dd{P}(K = k) \nonumber\\
  &= \dd{E}(K) + \dd{E}(X_{0}) + \sum_{\ell=1}^{\infty} \dd{E}(X_{\ell}) \dd{P}(K \ge \ell+1) \nonumber\\
  &=\frac 1q +  \dd{E}(X_{0}) + \frac {\lambda b p} {(1 - r)q} + \frac {(1-r) \dd{E}(X_{1}) - \lambda b } {(1 - r)(1 - p r)} p.
\end{align}

Under the scenario (2a), $\dd{E}(X_{k})$ of the $M/GI/1$ feedback queue will be used for the tail asymptotic of the sojourn time in the $GI/GI/1$ feedback queue. Thus, we introduce notations for them. Let $X_{k}^{(0)}(M/GI/1)$ be the $X_{k}$ of the $M/GI/1$ feedback queue for $X_{0} = 0$, then define $m^{(0)}_{k}$ as
\begin{align*}
  m^{(0)}_{0} = 0, \qquad m^{(0)}_{k} = \dd{E}(X_{k}^{(0)}(M/GI/1)), \quad k \ge 1.
\end{align*}
From \eq{X k mean}, we have
\begin{align}
\label{eq:m k}
 & m^{(0)}_{k} = \frac {1 - r^{k}} {1 - r} \lambda b = (1+r+\ldots + r^{k-1}) \lambda b .
\end{align}
We will use $m^{(0)}_{K-1}$ and $m^{(0)}_{K}$ in main results the next section. It should be noticed that they are random variables obtained by substituting $K-1$ and $K$ into $k$ of $m^{(0)}_{k}$.

\subsection{Main results}
\label{sect:main results}

We are ready to present the main results of this paper. They are proved in \sectn{proofs}.

\begin{theorem}
\label{thr:U asym 1}
For the stable $GI/GI/1$ feedback queue, assume that its service time distribution is intermediate regularly varying ($\sr{IRV}$). If the tagged customer finds the system empty, then
\begin{align}
\label{eq:U asym 1}
 \dd{P}\left(U > x\right) & \sim \frac 1q \dd{E}\left( 1 + X^{(0)}_{K-1} \right) \dd{P}\left((1+ m^{(0)}_{K-1}) \sigma > x \right)  \quad \mbox{ as } x \to \infty,
\end{align}
where $X^{(0)}_{k}$ is $X_{k}$ for $X_{0} = 0$, $m^{(0)}_{k} = \frac {1-r^{k}}{1-r} \lambda b$ for $r = p + \lambda b$ by \eq{m k}. Here the random
variable $K$ does not depend on $\{X_k^{(0)}\}$ and $\sigma$.
\end{theorem}

\begin{remark}
\label{rem:U asym 1}
In the case of a Poisson input, $\dd{E}\big(1+X^{(0)}_{K-1}\big) = \frac {q(1+p)} {1-rp}$. In general, it is hard to evaluate $\dd{E}\big(1+X^{(0)}_{K-1}\big)$ because this requires computing $\dd{E}(N^{e}(T_{k}))$.
\end{remark}

We prove this theorem in \sectn{proof U asym} using the PSBJ that is  established in \thr{GG1final}.

\begin{corollary}
\label{cor:U finite 1}
Under the assumptions of \thr{U asym 1}, for each $k \ge 1$,
\begin{align}
\label{eq:U finite 1}
  \dd{P}\left(T_{k} > x\right) & \sim \sum_{\ell=0}^{k-1} \dd{E}\left( 1 + X^{(0)}_{k-\ell-1} \right) \dd{P}\left((1+ m^{(0)}_{\ell}) \sigma > x \right) \quad \mbox{ as } x \to \infty.
\end{align}
\end{corollary}

This corollary is easily obtained from arguments used in the proof of \thr{U asym 1}. On the other hand, if we take the geometrically weighted sum of \eq{U finite 1} and if the interchange of this sum and the asymptotic limit are allowed, then we have \eq{U asym 1}. This interchange of the limits is legitimated by \thr{GG1final}. However, \cor{U finite 1} itself can be directly proved. We provide such a proof for a slightly extended version of \cor{U finite 1} in \app{alternative}. 

\begin{corollary}
\label{cor:M/GI/1}
Assume that the assumptions of \thr{U asym 1} hold. \\
(a) If the distribution of service times is $\sr{RV}$,
${\overline G}(x) = L(x)/x^{\alpha+1}$ with $\alpha >0$ and slowly varying function $L(x)$, then
$$
\dd{P}\left((1+ m^{(0)}_{K-1}) \sigma > x \right) \sim \frac{qL(x)}{x^{\alpha+1}}
\sum_{k=0}^{\infty} p^k \left(\frac {1-(p+ r^{k}\lambda b)} {1-r} \right)^{\alpha+1},
$$
where $r= p+\lambda b$.\\
(b) Under the assumption in (a), if the input stream is Poisson with parameter $\lambda$, then
$$
\dd{P} (U>x) \sim CL(x)x^{-(\alpha+1)}
$$ 
where
$$
C= \frac{q(1+p)}{(1-r)^{\alpha+1}(1-rp)} \sum_{k=0}^{\infty}
p^k (1-(p+r^k\lambda b))^{\alpha+1}.
$$
\end{corollary}

We next present the tail asymptotic for a tagged customer that arrives in the stationary system. By ``stationary'' we mean stationary in discrete time, i.e. at embedded arrival epochs, this is detailed in \sectn{psbj-st}.

\begin{theorem}
\label{thr:stationary case}

Let $U^0$ be the sojourn time of a typical customer in the stationary $GI/GI/1$ feedback queue with $\sr{IRV}$ distribution
$G$ of service times with mean $b$, $i.i.d.$ inter-arrival times with mean $a$ and probability of feedback $p=1-q \in (0,1)$. Let $\sigma_{I}$ be a random variable having the distribution function $G_{I}(x) \equiv 1 - \min\left(1,\int_{x}^{\infty} \ol{G}(u) du\right)$.\\
(a) Then, as $x\to\infty$,  
\begin{align}
\label{eq:final2-stat 2}
\dd{P}(U^0>x) & \sim \frac{\lambda}{1 - \lambda b} \left( \dd{P}(m^{(0)}_{K} \sigma_{I} > x) + \frac {(1-q) \rho} {1 - \rho} \dd{P}((\lambda b\, m^{(0)}_{K}) \sigma_{I} > x) \right),
\end{align}
where $m^{(0)}_{k}$ is defined by \eq{m k}, and $\sigma_{I}$ is independent of $K$.\\
(b) In particular, if $G$ is an $\sr{RV}$ distribution, $\overline{G}(x)=L(x)/x^{\alpha +1}$
with $\alpha >0$, where $L(x)$ is a slowly varying function, then
\begin{align}
\label{eq:final3-stat}
\dd{P}(U^0>x) & \sim \frac{L(x)}{x^{\alpha}} \cdot \frac{q(1-r)^{{-\alpha}}}{\alpha (a-b)}\left(
1+\frac{b(1-q)}{aq-b} \left(\frac{b}{a}\right)^{\alpha}\right)
\sum_{k=1}^{\infty} p^{k-1} (1-r^k)^{{\alpha}}. 
\end{align}
\end{theorem} 

\begin{remark}
\label{rem:stationary case}
(a) The tail function $\ol{G_{I}}(x) \equiv 1 - G_{I}(x)$ is called the ``integrated tail'' of $G$. Instead of $G_{I}$, we may use the stationary excess distribution $G_{e}(x) \equiv \frac 1b \int_{0}^{x} \ol{G}(u) du$ since $\ol{G}_{I}(x) = b \, \ol{G_{e}}(x)$ for sufficiently large $x$. In this case, let $\sigma_{e}$ be a random variable subject to $G_{e}$, then we can replace $\sigma_{I}$ by $\sigma_{e}$ in \eq{final2-stat 2}, multiplying its right-hand side by $b$.\\
(b) It is notable that no $X^{(0)}_{K-1}$ for the renewal arrivals is involved in \eq{final2-stat 2}. This is different from \eq{U asym 1}, and may come from averaging in the steady state.\\
(c) It may be interesting to compare the asymptotics in \eq{final2-stat 2} with those without feedback, which is well known (e.g., see \cite{BaFo2004}). Namely, let the stationary sojourn time $\widetilde{U}^0$ in the standard $GI/GI/1$ queue with inter-arrival times $\{t_n\}$ and with service times $\{\sigma_n^{H}\}$. where $\sigma_n^{H}$ has the same distribution as $\sum_{i=1}^{K} \sigma_{i}$. If $\sigma_{I}$ has a subexponential distribution, then
\begin{align}
\label{eq:final2-stat 3}
\dd{P}(\widetilde{U}^0>x) \sim \frac{\lambda}{1-\rho} \int_{x}^{\infty} \dd{P}\Big(\sum_{i=1}^{K} \sigma_{i} > u\Big) du \sim \frac{\lambda}{q(1-\rho)} \dd{P}(\sigma_{I} > x),
\end{align}
where the second asymptotic equivalence follows from Lh\'{o}pital's theorem and \eqref{Kesten}. Thus, one can see that \eq{final2-stat 2} is asymptotically compatible with \eq{final2-stat 3} as $q \uparrow 1$, because $\rho \to \lambda b$ and $m^{(0)}_{K} \to 1$ almost surely as $q \uparrow 1$.
\end{remark}

\section{Busy period and the principle of a single big jump}
\label{sect:psbj}
\setnewcounter

In this section, we present the Principle of a Single Big Jump (PSJB) in  \thr{GG1final} below, which will be used for a proof \thr{U asym 1}. For that, we first provide an auxiliary result on the tail asymptotics of the busy period in the $GI/GI/1$ queue without feedback. Denote its service time distribution by $H$ and let $\sigma^{H}_{i}$ be the $i$th service time. It is assumed that the arrivals are subject to the renewal process $N^{e}$ with interarrival times $t_i$ with mean $a$, and $H$ has a finite and positive mean $b_{H} > 0$. Denote the traffic intensity by $\rho \equiv b_{H}/a < 1$. Let $B$ be the (duration of the) first busy period in this $GI/GI/1$ queue, which is the time from the instant when the system becomes non-empty to the instant when it again becomes empty. We here omit the subscript $_{H}$ for $\rho, B$, because they will be unchanged for the $GI/GI/1$ feedback queue. We finally let $\tau^{H}$ be the number of customers served in the first busy period. 

We let $\xi^{H}_i=\sigma^{H}_{i}-t_{i}$, and let $S^{H}_{n} = \sum_1^n \xi^{H}_i$. Then $\tau^{H} = \min \{n\ge 1 : S^{H}_{n} \le 0\}$. Recall the definitions of classes
of heavy-tailed distributions ${\cal L}, {\cal S^*}, {\cal IRV}$ and ${\cal RV}$ at the end of \sectn{introduction}. The following theorem is proved in \sectn{proof gg1}.

\begin{theorem}
\label{thr:gg1}
Consider a stable $GI/GI/1$ queue, $\rho <1$,  with the service time distribution ${H}$.\\ 
If $H\in {\cal L}$, then
\begin{align}
\label{eq:B1}
\liminf_{x\to\infty}
\frac{\dd{P} (B>x)}{\dd{E} \tau^{H} \overline{H}(x(1-\rho ))}
\geq 1 \quad \mbox{and}
\quad \liminf_{x\to\infty}\frac{\dd{P} (\tau^{H} >x)}{\dd{E} \tau^{H} \overline{H}(x(a-b_{H}))}
\geq 1.
\end{align}
If, in addition, $H \in {\cal S}^{*}$, then,
for any $0<c<1$, 
\begin{align}
\label{eq:B2}
\limsup_{x\to\infty}
\frac{\dd{P} (B>x)}{\dd{E} \tau^{H} \overline{H}(c x(1-\rho ))}
\leq 1 \quad \mbox{and}
\quad \limsup_{x\to\infty}\frac{\dd{P} (\tau^{H} >x)}{\dd{E} \tau^{H} \overline{H}(c x(a-b_{H}))}
\leq 1.
\end{align}
Finally, if $H \in {\cal IRV}$, then, as $x\to\infty$,
\begin{align}
\label{eq:B3}
\dd{P} (B>x) \sim \dd{E} (\tau^{H}) \overline{H}(x(1-\rho ))
\quad \mbox{and}
\quad
\dd{P} (\tau^{H} >x) \sim \dd{E} (\tau^{H}) \overline{H}(x(a-b_{H})).
\end{align}
\end{theorem}

\begin{remark}
For the class of regularly varying tails, the equivalence \eq{B3}
was proved by Zwart in \cite{Zwar2001}. We provide a different proof which
is shorter and works for a broader class of distributions. Our proof
is based on probabilistic intuition related to the principle of a single big jump. A similar result holds for another class of distributions that overlaps
with the ${\cal IRV}$ class but does not contain it, see 
e.g. \cite{JeMoZw2004}.
\end{remark}

Recall the equivalence $A_{x} \simeq B_{x}$ for two families of events $A_{x}$ and $B_{x}$ with variable $x$. We have the following corollary, which is proved in \app{gg2}.

\begin{corollary}
\label{cor:gg2}
Consider the same $GI/GI/1$ queue as in \thr{gg1}. Let $\sigma^{H}_{n}$ be the service time of the $n$th arriving customer.
If $H\in {\cal IRV}$, then, for the busy period $B$, as $x \to \infty$,
\begin{align}
\label{eq:PSBJ-B1}
 & \dd{P} (B>x) \sim \sum_{n\ge 1} \dd{P} (\tau^{H}\ge n) \dd{P}(\sigma^{H}_n>x(1-\rho)) = \dd{E}\tau^{H} \dd{P}(\sigma^{H}_1>x(1-\rho)),
\end{align}
and, for any $\varepsilon > 0$, one can choose $N = N^{e}(\varepsilon) \ge 1$ such that, as $x \to \infty$,
\begin{align}
\label{eq:PSBJ-B2}
 & \dd{P} (B>x) \gtrsim \sum_{n=1}^N\dd{P} (\tau^{H} \ge n, \sigma^{H}_{n}\ge x(1-\rho )) 
\gtrsim (1-\varepsilon ) \dd{P} (B>x).
\end{align}
Furthermore, the following PSBJ holds:
\begin{equation}
\label{eq:PSBJ-B3}
\{B>x\} 
\simeq 
\cup_{n\ge 1} \{\tau^{H}\ge n, \sigma^{H}_{n} > x(1-\rho ) \}, \qquad x \to \infty.
\end{equation}
\end{corollary}

We now return to the $GI/GI/1$ feedback queue with the service time distribution $G$. Assume that the first customer arrives at the system at time instant $T_0=0$ and finds it empty. 
Recall that $K_i$ is the number of services $i$th customer has in the system, $K_{i}$'s are independent of everything else and $i.i.d$ with the same geometric distribution as $K$ (see \eq{K 1}). For convenience, we let $K=K_1$. Denote by $\sigma_{i}^{(j)}$ the $j$th service time of the $i$th customer in the $GI/GI/1$ feedback queue. Recall that $\sigma_{i}^{(j)}$ has the same distribution $G$.

Consider the $GI/GI/1$ queue without feedback and with service times $\sigma^{H}_i$ where
\begin{align*}
 \sigma^{H}_{i} = \sum_{j=1}^{K_{i}} \sigma_{i}^{(j)},
\end{align*}
and denote its distribution by $H$. Since the length of the busy period, $B$,
does not depend on the order of services, we may allow the server to proceed with services of lengths $\sigma_i^{j}$, like in the queue with feedback, and conclude that the (the lengths of) the busy periods are the same in both queues. Similarly, the traffic intensity $\rho$ in the new queue without feedback coincides with that in the $GI/GI/1$ queue with feedback. Furthermore, let $\tau$ be the number of service times in the first busy period of this feedback queue. Then, $\tau = \sum_{i=1}^{\tau^{H}} K_{i}$, and therefore we have
\begin{align*}
  \dd{E}(\tau) = \dd{E}(K) \dd{E}(\tau^{H}), \qquad b_{H} = \dd{E}(K) b.
\end{align*}

We now consider the $GI/GI/1$ feedback queue introduced in \sectn{GI/GI/1}.
We establish the PSBJ, i.e. show that, for large $x$, the rare event $\{U>x\}$ occurs mostly due to a big value of one of the service times.  
Our proof of Theorem \ref{thr:GG1final}  is based on \thr{gg1} and is given in \sectn{psbj 1}.

\begin{theorem}
\label{thr:GG1final}
Consider a stable single-server queue $GI/GI/1$ with feedback. Assume that the service times distribution is intermediate regularly varying. Denote by $U$ be the sojourn time of the first customer, and let
\begin{align}
\label{eq:P(x) 1}
  P_{k,\ell,i,j}(x) = \dd{P} (U>x, K=k+\ell, X_{k-1}=j,\sigma_{k,i}>x(1-\rho )).
\end{align}
If there exists a collection of positive functions $\{g_{k,\ell,i,j}(x)\}$ such that, as $x \to \infty$,
\begin{align}
\label{eq:GG1final g}
  P_{k,\ell,i,j}(x) \sim g_{k,\ell,i,j}(x), \qquad \forall k \ge 1, \ell \ge 0, 0\le i \le j,
\end{align}
and constants $C_{k,\ell,i,j}$ such that, for any $k \ge 1,\ell \ge 0, j \ge 0, 0 \le i \le j, x \ge 0$,
\begin{align}
\label{eq:g C}
 & g_{k,\ell,i,j}(x) \le C_{k,\ell,i,j} \cdot \dd{P}(\sigma > x),\\
\label{eq:C finite}
 & C:= \sum_{k=1}^{\infty}\sum_{\ell=0}^{\infty}\sum_{j=0}^{\infty}\sum_{i=0}^{j} C_{k,\ell,i,j} < \infty,
\end{align}
then 
\begin{equation}
\label{eq:GG1final2}
\dd{P} (U>x) \sim \sum_{k=1}^{\infty}\sum_{\ell=0}^{\infty}\sum_{j=0}^{\infty}\sum_{i=0}^{j} g_{k,\ell,i,j}(x).
\end{equation}
\end{theorem}

\section{Proofs of the theorems}
\label{sect:proofs}
\setnewcounter 

\subsection{Proof of \thr{gg1}}
\label{sect:proof gg1}

We will prove \thr{gg1} for the tail asymptotics of the busy period $B$ only. 
The proof for $\tau^{H}$, the number of arriving customers in the busy period, is similar. 

It is enough to prove the lower and upper bounds in \eq{B1} and \eq{B2}. Then the equivalences in \eq{B3} 
follow by letting $c$ tend to 1 and using the property of ${\cal IRV}$
distributions.

\noindent {\bf Lower bound.} 
Since $\xi^{H}_{1} = \sigma^{H}_{1} - t_{1}$, $d_0 \equiv b^{H}/(a-b^{H})$ is the solution to the equation 
\begin{align*}
  \dd{E} \sigma^{H}_1 + d_0 \dd{E} {\xi}^{H}_1 =0.
\end{align*}
We put $\widehat{\psi}^{H}_n = \sigma^{H}_n + d_{1} \xi^{H}_n$ for any positive number $d_{1} < d_{0}$. 
Then $\{\widehat{\psi}^{H}_n\}$ are i.i.d. random variables with common mean $\dd{E} \widehat{\psi}^{H}_1 >0$. Recall that $\tau^{H} = \min \{ n \ge 1 : S^{H}_n \leq 0\}$.
For any fixed real $C>0$, $R>0$, integer $N\ge 1$, and for $x\ge 0$, define events $D_{i}$ and $A_{i}$ for $i \ge 1$ as
\begin{align*}
 & D_i = \left\{ \sum_{j=1}^{i-1} |\widehat{\psi}^{H}_j|\le C, \tau^{H} \ge i, \widehat{\psi}^{H}_i>x+C+R \right\},
 \qquad A_i = \bigcap_{\ell \ge 1} \left\{\sum_{j=1}^{\ell} \widehat{\psi}^{H}_{i+j}\ge -R \right\}. 
\end{align*}
Then, we have
 \begin{equation}
 \label{double}
 \dd{P} (B>x) \ge 
\dd{P} \left(\sum_{i=1}^{\tau^{H}} \widehat{\psi}^{H}_i >x\right)  \ge \sum_{i=1}^N \dd{P} \left(D_i \cap A_i\right).
 \end{equation}
Here, the first inequality in \eqref{double} holds since $S^{H}_{\tau^{H}}$ is non-positive, and the second inequality comes from the following facts. Events $D_i$
are disjoint and, given the event $D_i$, we have $\sum_{j=1}^i \widehat{\psi}^{H}_j >x+R$. Then, given the event $D_i\cap A_i$,
we have $\sum_{j=1}^k \widehat{\psi}^{H}_j \ge x$ for all $k\ge i$ and, in particular, $\sum_{j=1}^{\tau^{H}} \widehat{\psi}^{H}_j>x$. Thus, \eqref{double} holds.

The events $\{A_i\}$ form a stationary sequence. Due to the SLLN, for any $\varepsilon >0$,  
one can choose $R$ so large that  $\dd{P}(A_i) \ge 1-\varepsilon$. For this $\varepsilon$ and any  $N \ge 1$, we can choose sufficiently large $C$ such that 
\begin{align*}
  \dd{P} \left(\sum_{j=1}^{i-1} |\widehat{\psi}^{H}_j|\le C, \tau^{H} \ge i\right) \ge (1-\varepsilon) \dd{P} (\tau^{H} \ge i), \qquad i \le N.
\end{align*}
Hence, \eqref{double} implies that, as $x\to\infty$,
\begin{align*}
\dd{P} (B>x) &\ge
\sum_{i=1}^N \dd{P} \left(\sum_{j=1}^{i-1} |\widehat{\psi}^{H}_j|\le C, \tau^{H} \ge i\right)
\dd{P}(\widehat{\psi}^{H}_i>x+C+R) \dd{P}(A_i) \\
&\ge
(1-\varepsilon )^2 \dd{P} (\widehat{\psi}^{H}_1 >x+C+R) \sum_{i=1}^N \dd{P} (\tau^{H} \ge i),
\end{align*}
and therefore the long-tailedness of distribution $H$ and (iii) of \rem{tailprop} yield
\begin{align*}
  \liminf_{x \to \infty} \frac {\dd{P} (B>x)}{\dd{P} (\widehat{\psi}^{H}_1 >x)} = \liminf_{x \to \infty} \frac {\dd{P} (B>x)}{\ol{H}(x/(1+d_{1}))} \ge (1-\varepsilon )^2 \sum_{i=1}^N \dd{P} (\tau^{H} \ge i).
\end{align*}
Letting first $N$ to infinity and then $\varepsilon$ to zero completes the proof of the first inequality of \eq{B1}.  

\noindent {\bf Upper bound.}
Take $L>0$ and put $\tilde{t}_n = \min (t_n,L)$,
$\tilde{\xi}^{H}_n = \sigma_n-\tilde{t}_n$,
$\tilde{S}^{H}_n = \sum_1^n \tilde{\xi}^{H}_i$,
$\tilde{\tau}^{H} = \min \{ n \ge 1: \tilde{S}^{H}_n \leq 0\}$.
Put also $\psi^{H}_n = \sigma^{H}_n + d_{2}\tilde{\xi}^{H}_n$ where $d_{2}>d_0$ is any number. 
Note that $\tilde{\xi}^{H}_1$ converges to $\xi^{H}_1$ in distribution and in mean, as $L\to\infty$.
We may choose $L$ so large that both $\dd{E} \tilde{\xi}^{H}_1$ and
$\dd{E} \psi^{H}_1$ are negative. Then $\tilde{\tau}^{H}$ is finite and
$\tilde{S}^{H}_{\tilde{\tau}^{H}} \in (-L,0]$ a.s.  
Further, as $L$ grows, $\tilde{\tau}^{H}$ converges to $\tau$ in distribution and 
in mean and, for any $\varepsilon >0$, we may choose $L$ so large that 
$\dd{E} \tau^{H} \le \dd{E}\tilde{\tau}^{H} \le (1+\varepsilon ) \dd{E} \tau^{H}$. 
By \rem{tailprop}, the distribution of $\psi^{H}_1$ also belongs to the class ${\cal S}^*$.
We have
\begin{align*}
  \dd{P} (B>x) &\leq 
\dd{P} \left(\sum_{i=1}^{\tilde{\tau}^{H}} \sigma^{H}_i > x\right) =
\dd{P} \left(\sum_{i=1}^{\tilde{\tau}^{H}} \psi^{H}_i > x +
d_{2}\tilde{S}^{H}_{\tilde{\tau}^{H}} \right)\\
&\leq 
\dd{P} \left(\sum_{i=1}^{\tilde{\tau}^{H}} \psi^{H}_i > x -d_{2}L\right) \le  \dd{P} \left(\max_{1\le j\le \tilde{\tau}^{H}} \sum_{i=1}^j \psi^{H}_i > x-d_{2}L \right)\\
&\sim \dd{E} \tilde{\tau}^{H} \dd{P} (\psi^{H}_1 > x-d_{2}L)
\leq (1+\varepsilon )\dd{E} \tau^{H} \dd{P} (\psi^{H}_1 > x-d_{2}L),
\end{align*}
where the equivalence follows from \thr{FZ1}. Further, 
$$
 \dd{P} (\psi^{H}_1 > x-d_{2}L) 
\sim 
\dd{P} (\psi^{H}_1 > x) \\
\sim 
\dd{P} ((d_{2}+1)\sigma^{H}_1 > x) =
\overline{H}(x/(d_{2}+1))
$$
where the first equivalence follows from the long-tailedness of the distribution of $\psi^{H}_1$ and the second from \rem{tailprop}. 
Letting $\varepsilon$ tend to zero, we have 
$$
\lim \sup_{x\to\infty}
\frac{\dd{P} (B>x)}{\dd{E} \tau^{H} \overline{H}(x/(d_{2}+1))}
\leq 1
$$
for any $d_{2} > d_0$. Let $c = (d_0+1)/(d_{2}+1) <1$. This proves the first inequality in \eq{B2}.

\subsection{Proof of \thr{GG1final}}
\label{sect:psbj 1} 

Recall that we consider the  scenario where the initial customer $1$  
arrives at the empty system. Clearly, $\sigma^{H}_1 \le U \le B$ a.s. where 
$\sigma^{H}_1$
is the total service time of customer $1$, and $B$ is the duration of the first busy period. 

Equivalence relation \eqref{Kesten} and \rem{tailprop} 
from the Appendix 
imply that, given $\sigma$ has an intermediate regularly varying distribution, random variables $\sigma^{H}_i$ have 
a common intermediate varying distribution too.  Since
any intermediate varying distribution is dominantly varying (see Property \eq{class order}), 
we get from \thr{gg1} that  
\begin{equation}\label{weakTE}
\limsup_{x\to\infty} \dd{P}(B>x)/\dd{P}(\sigma^{H}_1 >x) <\infty.
\end{equation}
Relations  
\eqref{weakTE} and \eqref{Kesten} lead then to 
the logarithmic asymptotics: 
\begin{lemma}
 Assume that the distribution of the service time $\sigma_{1}^{(1)}$ belongs to the class ${\cal IRV}$. Then,
as $x\to\infty$,
\begin{equation}\label{heavy-log}
\log \dd{P} (U>x) \sim \log \dd{P} (\sigma^{H}_1>x) \sim
\log \dd{P} (\sigma_{1,1}>x).
\end{equation}
\end{lemma}

Further, by \cor{gg2}, the PSBJ for $B$ holds:
\begin{equation}\label{BP}
\dd{P} (B>x) \sim 
\dd{P} \left(B >x, \cup_{i=1}^{\tau^{H}} \{\sigma^{H}_i >x (1-\rho )\}\right)
\sim \dd{P} \left(\cup_{i=1}^{\tau^{H}} \{\sigma^{H}_i >x (1-\rho )\}\right)
\end{equation}
Here $\tau$ is the number of customers served within the first busy period.

Combining \eqref{BP} and \eqref{Kesten}, we arrive at the following result:

\begin{lemma}\label{PSBJJ}
Consider a stable single-server queue $GI/GI/1$ with feedback. Let $B$ be the duration of the first busy period and 
$U$ the sojourn time of the first customer. Assume that the service times distribution is intermediate regularly varying.
Then 
\begin{equation}\label{IRVEA1}
\dd{P} (U>x) = \dd{P} (U>x,B>x)\sim \dd{P} \left(\{U>x\} 
\bigcap  \bigcup_{i=1}^{\tau^{H}}
 \bigcup_{j=1}^{K_i}\{\sigma_{i}^{(j)}>x (1-\rho)\}\right).
\end{equation}
\end{lemma}

To derive the exact asymptotics for $\dd{P} (U>x)$, we recall that, for $1\le k <K \equiv K_{1}$, $X_{k} \ge 0$ is the total number of services of other customers between the $k$th and the $(k+1)$st services
of customer 1, and let $\sigma_{k,i}$ be the service time of the $i$ service there, $1\le i \le X_{k}$. Further, under the scenario (\sect{main}a), $X_0=0$.  Then 
let $\nu \ge 0$ be the total number of services of other customers after the departure of the first customer within the
busy period, and let $\sigma_{i}^{*}$ be the $i$th service time there, $1\le i \le \nu$. 
Then
random variables $\sigma_{k,.i}$ and $\sigma_{i}^{*}$ are $i.i.d.$ with the same distribution as $\sigma$ and $U$ is given by \eq{U 1 (a)}. From \eqref{IRVEA1}, we get,
\begin{align*}
\dd{P} (U>x) \sim \dd{P} \left(\{U>x\} 
\bigcap  \left( \bigcup_{i=1}^{K_1}\bigcup_{j=1}^{X_i}
\{\sigma_{j}^{(i)}>x (1-\rho)\}\bigcup\bigcup_{j=1}^{\nu} \{\sigma_j^{*}>x(1-\rho )\}\right)\right).
\end{align*}
On the other hand, we have 
\begin{align*}
\dd{P} \left(\{U>x\} 
\bigcap  \bigcup_{i=1}^{\nu} \{\sigma_i^{*}>x(1-\rho )\}\right) & =
\sum_{n\ge 1} \dd{P} (U>x,\nu =n) \cdot \dd{P}\left(
\bigcup_{i=1}^{n} \{\sigma_i^{*}>x(1-\rho )\}\right)\\
&\le  C \sum_{n\ge 1} \dd{P} (U>x,\nu =n) n \overline{G}(x)\\
&= C \dd{E} \left( \nu \cdot {1} (U>x) \right) 
\overline{G}(x)\\
&\le   C \dd{E} \left(\tau \cdot {1} (U>x)\right)
\overline{G}(x) = o(\overline{G} (x)),
\end{align*}
where $C=\sup_x \overline{G}(x(1-\rho))/\overline{G}(x)<\infty$ is a  constant, and recall that $\tau$ is the total number of services within the busy cycle. The last line follows since 
$\dd{E} \tau = \dd{E} \tau^{H} /q$ is finite. 
Therefore, we have
\begin{lemma}\label{PSBJJ2}
In the conditions of Lemma \ref{PSBJJ}, we have
\begin{equation}\label{IRVEA2}
\dd{P} (U>x) \sim \dd{P} \left(\{U>x\} 
\bigcap  \bigcup_{k =1}^{K_1} \bigcup_{i=0}^{X_{k-1}}
\{\sigma_{k,i}>x (1-\rho)\}\right).
\end{equation}
\end{lemma}

Moreover, the following result holds:
\begin{lemma}\label{PSBJJ3}
Assume the conditions of Lemma \ref{PSBJJ} hold. 
Then, for any $\varepsilon >0$, one can find $N$ such that
\begin{equation}\label{IRVEA3}
\dd{P} (U>x)  \gtrsim \dd{P} (D_{N}(x))\gtrsim (1-\varepsilon ) \dd{P} (U>x)
\end{equation}
where 
\begin{equation}\label{DNL}
D_{N}(x) = 
\bigcup_{k=1}^N
\bigcup_{\ell=0}^{2N}
\left\{\{U>x, K_1=k+\ell\} 
\bigcap \bigcup_{i=0}^{X_{k-1}}
\{\sigma_{k,i}>x (1-\rho)\} \right\}.
\end{equation}
\end{lemma}
{\sc Proof.} Indeed, the term in the right-hand side of \eqref{IRVEA2} is bigger
than $\dd{P} (D_{N}(x))$ and smaller than the sum 
$\dd{P} (D_{N}(x)) + \dd{P} (U>x, K_1 >N)$, where 
\begin{equation}\label{IRVEA4}
\dd{P} (U>x, K_1 >N) \le \dd{P} (B>x, K_1>N).
\end{equation}
Consider again the auxiliary $GI/GI/1$ queue with service times $\sigma^{H}_i
=\sum_{j=1}^{K_i}\sigma_{i}^{(j)}$ and the first-come-first-served service discipline.
Consider the following majorant: assume that at the beginning of the first cycle, in addition to customer 1, an extra $K-1$ new customers arrive, so there are $K$ arrivals in total. Here $K$ is a geometric random variable with parameter $p$ that does not depend on service times.
Then the first busy period in this queue has the same
distribution as $\sum_{i=1}^K B_i$ where $B_i$ are $i.i.d.$ random variables that have the same distribution as $B$ and do not depend on $K$.
By monotonicity,
$$
\dd{P} (B>x, K_1>N) \le \dd{P} \left(\sum_{i=1}^K B_i>x, K>N\right) 
=\dd{P} \left(\sum_{i=1}^{K{1}(K>N)} B_i >x\right).
$$
Due to \eqref{Kesten}, the latter probability is equivalent, as $x\to\infty$, to
$$
\dd{E} (K{1} (K>N)) \dd{P}(B>x)\le
C_0\dd{E} (K{1} (K>N))\dd{E} K \overline{G}(x)
$$
where $C_0$ is from \eqref{weakTE}. Now choose $N$ such that 
$C_0 \dd{E} (K{1} (K>N))\dd{E} K \le \varepsilon$.
Since $\dd{P}(U>x) \ge \overline{G}(x)$, \eqref{IRVEA3} follows.

We can go further and obtain the following result.
\begin{lemma}
\label{lem:PSBJJ4}
Assume the conditions of Lemma \ref{PSBJJ} hold, and let
\begin{equation}\label{GNLR}
G_{N,R}(x) = \bigcup_{k=1}^N \bigcup_{\ell=0}^{2N}
\left\{\{U>x, K_1=k+\ell\} \bigcap \bigcup_{i=0}^{R} \{\sigma_{k,i}>x (1-\rho)\} \right\},
\end{equation}
Then, for any $\varepsilon >0$, 
one can choose a positive integer $R$ such that
\begin{equation}
\label{IRVEA5}
\dd{P} (D_{N}(x)) 
\ge
\dd{P} (G_{N,R}(x)) 
\ge
(1-\varepsilon)\dd{P} (D_{N}(x))
\end{equation}
where the event $D_{N}(x)$ was defined in \eqref{DNL}.
Further, 
\begin{equation}
\label{IRVEA6}
 \dd{P} (G_{N,R}(x)) \sim
 \sum_{k=1}^N\sum_{\ell=0}^{2N}\sum_{j=0}^R\sum_{i=0}^{j}
 \dd{P} (U>x, K_1=k+\ell, X_{k-1}=j,\sigma_{k,i}>x(1-\rho )).
 \end{equation}
\end{lemma}
{\sc Proof}. 
Indeed,
\begin{align*}
\dd{P} & \left(D_{N}(x) \setminus G_{N,R}(x)\right)\\
&\le \sum_{k=1}^{N} \sum_{\ell=0}^{2N} \sum_{j > R}
 \dd{P} \left( \{U>x, K_1=k+\ell, X_{k-1} = j\} \bigcap \bigcup_{i=R+1}^{j} \{\sigma_{k,i}>x (1-\rho)\} \right) \\
 &\le 
 \sum_{k=1}^{N} \sum_{j>R} \sum_{i=0}^{j} 
 \dd{P} ( X_{k-1}=j, \sigma_{k,i}>x(1-\rho))\\
 &=
 \sum_{k=1}^{N}\sum_{j>R} (j+1) \dd{P} ( X_{k-1}=j) 
 \overline{G}(x(1-\rho ))\\
 &= \sum_{k=1}^{N}\dd{E} (( X_{k-1}+1) {1}( X_{k-1} >R))
 \overline{G}(x(1-\rho )) \\
 &\le 
 N \dd{E} ((\tau+1){1} (\tau>R)) \overline{G}(x(1-\rho ))
 \end{align*}
 where the term $\dd{E} ((\tau+1){1} (\tau>R))$ may be made 
 as small as possible by taking a sufficiently large $R$.  
 Then \eqref{IRVEA6} follows 
since the probability of a union of events is always smaller than
the sum of their probabilities, and is bigger than the sum of probabilities of events minus the sum of probabilities of pairwise
intersections of events. Each probability of intersection of two independent 
events is smaller than
$$
\dd{P} (\sigma_1^{(1)}>x(1-\rho),\sigma_2^{(1)}>x(1-\rho)) \le C\overline{G}^2(x) = o(\overline{G}(x)),
$$
therefore their finite sum is  $o(\overline{G}(x))$ and \eqref{IRVEA6}
follows.

We are now in a final step of the proof of \thr{GG1final}. For $k \ge 1,\ell, j \ge 0$, define $D_{k,\ell,j}$ as
\begin{align}
\label{eq:D klj}
  D_{k,\ell,j} = \dd{P} (K=k+\ell, X_{k-1}=j) = \dd{P} (K>k, X_{k-1}=j) \dd{P} (K=\ell),
\end{align}
where the second equality holds because $K$ is geometrically distributed. Then, \lem{PSBJJ4} implies \eq{GG1final2} for $g_{k,\ell,i,j}(x) = P_{k,\ell,i,j}(x)$
since, for any $k,\ell,j\ge i$, 
$$
g_{k,\ell,i,j}(x) \le \dd{P} (K=k+\ell, X_{k-1}=j, \sigma_{k,i}>x(1-\rho)) = D_{k,\ell,j} \dd{P}(\sigma >x(1-\rho)),
$$
and 
\begin{align}
\label{eq:D finite}
  \sum_{k,\ell,i \le j} D_{k,\ell,j} = \sum_{k,\ell,j} jD_{k,\ell,j} = \sum_{k,j} j \dd{P}(X_{k-1}=j, K>k) \le \dd{E} \tau/q <\infty,
\end{align}
where, recall, $\tau$ is the total number of customers served in the first busy period. Clearly, \eq{GG1final2} is also valid for a general $\{g_{k,\ell,i,j}(x)\}$ because of the conditions \eq{g C} and \eq{C finite}. This completes the proof of \thr{GG1final}.

\subsection{Proof of \thr{U asym 1}}
\label{sect:proof U asym}
 
We first recall the notation: $U_{1}, U_{2},\ldots$ and $X_{0}, X_{1}, \ldots$ are the service cycles and the number of customers other than the tagged customer served in the cycles, respectively. Here $u_{0} = X_{0} = 0$. 
In general, the sojourn time is a randomly stopped sum of $i.i.d.$ positive random variables, and both  the summands and the counting random variable have heavy-tailed distributions. It is known that 
it is hard to study the tail asymptotics for general heavy-tailed distributions (see, e.g., \cite{Gree1973}) in this case. We proceed under the assumption that the service time distribution is intermediate regularly varying.

Recall that $\sigma_{k,0}$ is the $k$th service time of the tagged customer and, for $i=1,\ldots, X_k$, $\sigma_{k,i}$ is the $i$th service time in the queue $X_k$. Further, $T_{k}=\sum_{\ell=1}^k U_{\ell}$ be the time instant when the $k$th service of the tagged customer is completed, where $U_{1} = \sigma_{1,0}$. Introduce the notation 
\begin{align*}
  U_k^+ = \sum_{\ell=k+1}^K U_\ell, \qquad k \ge 1,
\end{align*}
which is the remaining time the tagged customer spends in the system after the completion of the $k$th service, and let $v_{k}$ be the residual inter-arrival time of the input when the $k$th service of the tagged customer ends.

In what follows, we will say that an event involving some constants and functions/sequences 
occurs ``with high probability'' if, for any $\varepsilon >0$, there exists constants and functions/sequences (that depend on $\varepsilon$) with the desired properties such that the event occurs with probability at least $1-\varepsilon$.

For example, let $S^{\sigma}_n = \sum_1^n \sigma_i$ be the sum of $i.i.d.$ random variables with finite mean $b$. Then the phrase 
``with high probability (WHP), for all $n=1,2,\ldots$,
$$ 
S^{\sigma}_n \in (n(a-\delta_n)-C, n(a+\delta_n)+C)
$$
with $C>0$ and $\delta_n\downarrow 0$'' 
means that
``for any $\varepsilon >0$, there exist a constant $C\equiv C_{\varepsilon}>0$ and a sequence $\delta_n\equiv\delta_n(\varepsilon) \downarrow 0$ such that the probability of the event
$$
\{ S^{\sigma}_n \in (n(a-\delta_n)-C, n(a+\delta_n)+C), \ \mbox{for all} \ \ n\ge 1\}
$$
is at least $1-\varepsilon$''. We can say equivalently that 
``WHP, for all $n=1,2,\ldots$,
$
S^{\sigma}_n \in (na - o(n), na+ o(n))
$
'', or, simply, ``WHP, $S^{\sigma}_n \sim an$'', 
and this means that 
``for any $\varepsilon >0$, there exists a positive function $h(n)=h_{\varepsilon}(n)$ which is an $o(n)$-function 
(it may tend to infinity, but slower than $n$) and is such that
the probability of the event
$$
\{ S^{\sigma}_n \in (na - h(n), na+ h(n)), \ \mbox{for all} \ n\}
$$
is at least $1-\varepsilon$.''


Now we show \eq{GG1final g} for
\begin{align}
\label{eq:g 1}
  g_{k,\ell,i,j}(x) & = \dd{P}\left( K = k+\ell, X_{k-1} = j\right) \dd{P}((1+m^{(0)}_{\ell}) \sigma_{k,i} > x), \nonumber\\
  & \hspace{30ex} k\ge 1, \ell \ge 0, j\ge 0, 0\le i\le j.
\end{align}
Namely, we show that, for all $k \ge 1, \ell \ge 0, 0\le i \le j$,
\begin{align}
\label{eq:g 2}
  \dd{P} (U>x, K=k+\ell, X_{k-1}=j,\sigma_{k,i}>x(1-\rho )) \sim g_{k,\ell,i,j}(x).
\end{align}

We prove \eq{g 2} by induction on $\ell \ge 0$, for each fixed $k\ge 1, 0 \le i \le j$.

\noindent {\bf Lower bound, $\ell=0.$}\\
Since $\sigma_{k,i} > x$ implies that $U > x$ and $\sigma_{k,i} > (1-\rho)x$, the lower bound for the LHS of \eq{g 2} is
\begin{align*}
\dd{P} (K=k, X_{k-1}=j) \overline{G}(x) = g_{k,0,i,j}(x),
\end{align*}
because $m^{(0)}_{0} = 0$.

\noindent {\bf Upper bound, $\ell=0$.} \\
There is a constant $w > 0$ such that $T_{k-1}\le w$ and $\sum_{0\le i' \le j, i' \ne i} \sigma_{k,i'} \le w$ WHP.
Then $U\le 2w+\sigma_{k,i}$, so the upper bound  for the LHS of \eq{g 2} is
\begin{align*}
&
\varepsilon\overline{G}(x(1-\rho)) + (1+o(1)) \dd{P} (K=k, X_{k-1}=j, \sigma_{k,i}+2w>x)\\
& \quad \sim 
\varepsilon\overline{G}(x(1-\rho)) + \dd{P} (K=k, X_{k-1}=j) \overline{G}(x).
\end{align*}
Letting $\varepsilon$ tend to zero in this upper bounds yields that the lower and upper bounds are asymptotically identical. Since $m^{(0)}_{0} = 0$, they are further identical to $g_{k,0,i,j}(x)$ of \eq{g 1}. Thus, \eq{g 2} is verified.

Turn to the case $\ell=1.$

\noindent {\bf Lower bound, $\ell=1.$}\\ 
Like in the case $\ell=0$, replace all other service times $\sigma_{k,i'}, i'\ne i$
by zero. Assume that all $j$ customers from the group $X_{k-1}$ leave the system after their service completions. 
WHP, $v_{k-1}\le w$. Given $y=\sigma_{k,i}$ is large and much bigger than $w$, we have that at least
$N^{e}(y-w)$ customers arrive during time $U_{k} \ge \sigma_{k,i} = y$. Again WHP,
$$
N^{e}(y-w) \in (\lambda y - o(y), \lambda y + o(y))
$$
and, again WHP, their total service time is within the time interval 
$(\lambda b y - o(y), \lambda b y + o(y))$. Therefore,
$$
U\ge U_{k}+U_{k+1} \ge y + \lambda b y - o(y)
$$
and the RHS is bigger than $x$ if $y>x/(1+\lambda b) + o(x)$.
Therefore,  the lower bound for the LHS of \eq{g 2}
is
\begin{align*}
&(1+o(1)) \dd{P} (K=k+1, X_{k-1}=j, \sigma_{k-1,i}>x/(1+\lambda b) + o(x)) - \varepsilon\overline{G}(x(1-\rho))\\
& \quad \sim 
\dd{P} (K=k+1, X_{k-1}=j) \overline{G}(x/(1+\lambda b)) - \varepsilon\overline{G}(x(1-\rho)).
\end{align*}
 
\noindent {\bf Upper bound, $\ell=1.$} \\
WHP, $T_{k-1}\le w$, $v_k\le w$ and $\sum_{0 \le i' \le j, i'\ne i} \sigma_{k,i'} \le w$. Let $\sigma_{k,i}=y \gg 1$.
Then, WHP, $U_{k}\le y+2w$ and the number of external arrivals within
$U_{k}$ is bounded above by $1+N^{e}(y+2w) = \lambda y + o(y)$, again WHP. 
Assume that all $X_{k-1}=j$ customers stay in the system after their services. Then again
$j+1+N^{e}(y+w) = \lambda y +o(y)$, WHP. Therefore,
$U_{k+1} = b \lambda y + o(y).$ Then we arrive at the upper bound that meets the
lower bound. 

Thus, \eq{g 2} is verified for $g_{k,2,i,j}(x)$ because $m^{(0)}_{1} = \lambda b $ by \eq{m k}.

\noindent {\bf Induction step.}\\
We can provide induction for any finite number of steps. Here is the induction base.
Assume that $\sigma_{k,i}=y \gg 1$ and that, after $\ell \ge 1$ steps, $T_{k+\ell'} \sim (1+m^{(0)}_{\ell'}) y$ for $0 \le \ell' \le \ell$, and there are
$X_{k+\ell-1}$ customers in the queue and that $X_{k+\ell-1}=wy + o(y)$, WHP, 
where $w>0$. Then, combining upper and lower bounds, we may conclude that,
again WHP, 
$U_{k+\ell} = b wy + o(y)$ and then
\begin{align}
\label{eq:Tk 1}
 & T_{k+\ell} = T_{k+\ell-1} + U_{k+\ell} \sim (1+m^{(0)}_{\ell-1} + bw) y,\\
\label{eq:Xk 1}
 & X_{k+\ell} = pX_{k+\ell-1} + \lambda U_{k+\ell} + o(y) = wy r + o(y),
\end{align}
where we recall that $r = p + \lambda b$. By the induction hypothesis, \eq{Tk 1} implies that
\begin{align*}
  1+m^{(0)}_{\ell} = 1+m^{(0)}_{\ell-1} + bw,
\end{align*}
which, with \eq{m k}, yields that
\begin{align*}
  bw = m^{(0)}_{\ell} - m^{(0)}_{\ell-1} = r^{\ell-1} \lambda b.
\end{align*}
Hence, by \eq{Xk 1},
\begin{align*}
  T_{k+\ell+1} = T_{k+\ell} + U_{k+\ell+1} & \sim (1+m^{(0)}_{\ell}) y + bw ry \nonumber\\
  & = (1 + m^{(0)}_{\ell} + r^{\ell} \lambda b) y = (1 + m^{(0)}_{\ell+1}) y.
\end{align*}
This completes the induction step for $\ell+1$. 

We finally check the conditions \eq{g C} and \eq{C finite}. Since an intermediate regularly varying distribution is dominantly varying, it follows from \eq{g 1} that
\begin{align*}
  g_{k,\ell,i,j}(x) & \le \dd{P}\left( K = k+\ell, X_{k-1} = j\right) \dd{P}((1+m^{(0)}_{\infty}) \sigma  > x)\\
  & \le \dd{P}\left( K = k+\ell, X_{k-1} = j\right) c \dd{P}(\sigma  > x) 
\end{align*}
for some $c > 0$. Hence, letting
\begin{align*}
  C_{k,\ell,i,j} = \dd{P}\left( K = k+\ell, X_{k-1} = j\right) c,
\end{align*}
\eq{g C} is verified, while \eq{C finite} follows from \eq{D finite}. Thus, by \thr{GG1final},
\begin{align*}
 \dd{P}(U > x) & \sim \sum_{k=1}^{\infty}\sum_{\ell=0}^{\infty}\sum_{j=0}^{\infty}\sum_{i=0}^{j} \dd{P}\left( K = k+\ell, X_{k-1} = j\right) \dd{P}\left((1+m^{(0)}_{\ell}) \sigma > x \right)\\
 & = \frac 1q \sum_{k=1}^{\infty}\sum_{\ell=0}^{\infty} qp^{k}\dd{E}\left( 1 + X_{k-1}\right) qp^{\ell-1} \dd{P}\left((1+m^{(0)}_{\ell}) \sigma > x \right),
\end{align*}
which implies \eq{U asym 1}, and \thr{U asym 1} is proved.

\subsection{PSBJ for the stationary queue}
\label{sect:psbj-st}

We now consider the case where customer $1$ arrives to the stationary queue and denote by $U^0$ its sojourn time. 
In this Section, we frequently use the following notation: for a distribution $F$ having a finite mean,  $\overline{F_{I}}(x) = \min (1, \int_x^{\infty} \ol{F}(y) dy)$ is its integrated tail distribution (see (a) of \rem{stationary case}).

\begin{figure}[h] 
   \centering
   \includegraphics[height=4.5cm]{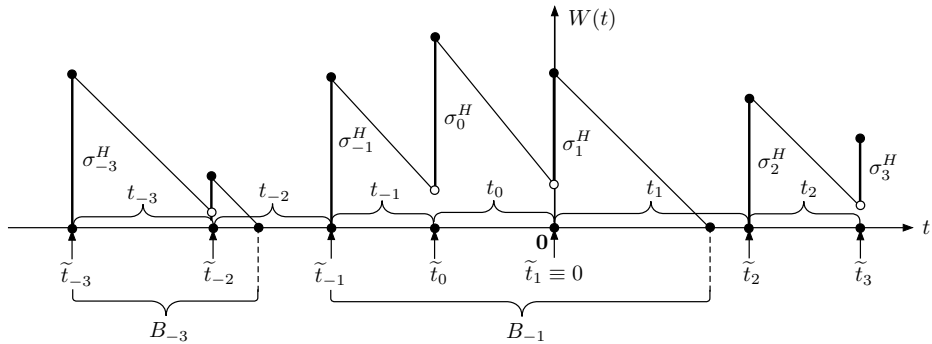}
\caption{Workload and time indexes under the stationary regime: $B^{0} = B_{-1}+\widetilde{t}_{-1}$ and $\arg_{k} \{\min_{k \le \ell  \le 0} S^{H}_{k,\ell }>0\} = -1$. At and after time $0$, the time indexes are identical with those in the previous sections.}
\label{fig:Feedback_time_indeces}
\end{figure}

By ``stationarity'' we mean stationarity in discrete time, i.e. at embedded arrival epochs.
So we assume that the system has started from time $-\infty$ and that customer $1$ arrives at time $\widetilde{t}_{1} \equiv 0$, customers with indices $k \le 0$ 
enter the system at time instants $\widetilde{t}_k=-\sum_{j=k}^0 t_j$ and customers with indices $k \ge 2$ 
at time instants $\widetilde{t}_k=\sum_{j=1}^{k-1} t_j$.
For $k \le \ell$, let $S^{H}_{k,\ell} = \sum_{j=k}^{\ell}\xi^{H}_j$, where $\xi^{H}_{j} = \sigma^{H}_{j} - t_{j}$. Then the stationary busy cycle covering  $0$ starts at $\widetilde{t}_k$, 
$k \le 0$ if 
$$
\{\sup_{j<k} S^{H}_{j,k}\le 0, \ \min_{k \le \ell  \le 0} S^{H}_{k,\ell }>0\}.
$$
So, if $B^0$ is the remaining duration of the busy period viewed at time $0$, then
\begin{align}\label{eq:decomp}
\{B^0>x\} & = \bigcup_{k = -\infty}^{0}  \left\{\sup_{j<k} S^{H}_{j,k}\le 0, B_k > -\widetilde{t}_k + x\right\},
\end{align}
where $B_k$ is the duration of the period that starts at time $\widetilde{t}_k$ given that  customer 
$k$ arrives in the empty system (then, in particular, $B=B_0$). See \fig{Feedback_time_indeces}.

Let 
\begin{align*}
  \tau^{H}_{k} = \min\left\{ n \ge 1; S^{H}_{k,n} \le 0 \right\}, \qquad k \le 0,
\end{align*}
which is the number of customers arriving at or after time $0$ in the busy period when it starts at time $\widetilde{t}_{k}$, and let
\begin{align*}
  A^{H}_{k} = \{\sup_{j <k} S^{H}_{j,k}\le 0\}, \qquad k \le 0,
\end{align*}
then $A^{H}_{k} \cap \{\tau^{H}_{k} \ge \ell\}$'s for $k \le 0, \ell \ge -k+1$ are disjoint sets.

%

Thus, \eq{decomp} may be written as
\begin{align}
\label{eq:decomp2}
  \{B^0>x\} & = \bigcup_{k = 0}^{\infty} 
  \left(A^{H}_{-k} \cap \left\{\tau^{H}_{-k} \ge k+1, B_{-k} > -\widetilde{t}_{-k} + x\right\}\right),
\end{align}
and therefore, applying the PSBJ of \cor{gg2} to each busy period $B_{-k}$,
\begin{align*}
 & \{B^0>x\} = \bigcup_{k = 0}^{\infty} \left(A^{H}_{-k} \cap \left\{\tau^{H}_{-k} \ge k+1, B_{-k} > -\widetilde{t}_{-k} + x\right\}\right) \nonumber\\
 & \quad \simeq \bigcup_{k = 0}^{\infty} \bigcup_{i = 1}^{\infty} \left(A^{H}_{-k} \cap \left\{\tau^{H}_{-k} \ge \max(k+1,i), \sigma^{H}_{-k+i} > (-\widetilde{t}_{-k} + x)(1-\rho) \right\}\right),
\end{align*}
where $\sigma^{H}_{-k+i}$, $i\ge 0$ is the service time of the $i$-th customer arriving in the busy period that starts at time $\widetilde{t}_{-k}$.

Hence, letting
\begin{align*}
 A^{0}_{-}(x) & = \bigcup_{k=1}^{\infty} \bigcup_{i=1}^{k} \left(A^{H}_{-k} \cap \left\{\tau^{H}_{-k} \ge k+1, \sigma^{H}_{-k+i} > (-\widetilde{t}_{-k}+x)(1-\rho)\right\} \right) \nonumber\\
 A^{0}_{+}(x) & = \bigcup_{k=0}^{\infty} \bigcup_{i=k+1}^{\infty} \left(A^{H}_{-k} \cap \left\{\tau^{H}_{-k} \ge i,\sigma^{H}_{-k+i} > (-\widetilde{t}_{-k}+x)(1-\rho)\right\} \right),
\end{align*}
we have
\begin{align}
\label{eq:decomp3}
  \{B^0>x\} & \simeq A^{0}_{-}(x) \cup A^{0}_{+}(x).
\end{align}

We first consider the event $A^{0}_{+}(x)$, which is a contribution of big jumps at or after time $0$, and show that its probability is negligible with respect to
$\overline{H_I}(x)$, as $x\to\infty$.
Clearly, for any positive function $h(x)$ and 
for any $\varepsilon\in (0,a)$, 
\begin{align*}
A_+^0(x) & \subseteq \bigcup_{k\ge 0} \{-\widetilde{t}_{-k} < (a-\varepsilon )k - h(x)\}\\
& \qquad \bigcup \bigcup_{k\ge 0} \bigcup_{i\ge k+1}  
\{\tau_{-k}^H\ge i, \sigma_{-k+i}^H > ((a-\varepsilon)k+x-h(x))(1-\rho))\}.
\end{align*}
Then
\begin{align*}
\dd{P}(A_+^0(x)) & \le \sum_{k=0}^{\infty} \sum_{i=k+1}^{\infty}
\dd{P}(\tau_{-k}^H\ge i, \sigma_{-k+i}^H > (-\widetilde{t}_{-k}+x)(1-\rho ))\\
& \le \sum_{k\ge 0} \dd{P}(-\widetilde{t}_{-k}<(a-\varepsilon )k - h(x)) \\
& \qquad + \sum_{k\ge 0} \sum_{i\ge k+1} \dd{P} (\tau_{-k}^H\ge i) \dd{P} (\sigma^H_{-k+i}>
((a-\varepsilon)k + x - h(x))(1-\rho))\\
& \le C e^{-\alpha h(x)} + \sum_{k\ge 0} \dd{E} ((\tau^H-k+1)^+) \dd{P}
 (\sigma^H_{-k+i}>
((a-\varepsilon)k + x - h(x))(1-\rho))\\
& = o(\overline{H_I}(x(1-\rho)))= o(\overline{G_I}(x(1-\rho))),
\end{align*}
if one takes, say, $h(x) = x^c$ for some $c<1$. Here the second inequality follows since
 $\left\{\tau^{H}_{-k+i} \ge i\right\} = \left\{\tau^{H}_{-k+i} \le i-1\right\}^{c}$ is independent of $\sigma^{H}_{-k+i}$, the 
third inequality from Chernoff's inequality, for a small $\alpha>0$, and the final conclusion from 
property \eqref{PSBJ101} in the Appendix.

Thus, we only need to evaluate the contribution of big jumps that occur before time $0$. Namely, we analyse $A^{0}_{-}(x)$.
Note that, for any $k_0>0$, the probability of the event 
\begin{align*}
A_-^{(0, k_0)}(x) & = \bigcup_{k=1}^{k_0} \bigcup_{i=1}^{k} \left(A^{H}_{-k} \bigcap \left\{\tau^{H}_{-k} \ge k+1, \sigma^{H}_{-k+i} > (-\widetilde{t}_{-k}+x)(1-\rho)\right\} \right)
\end{align*}
is of order $O(\overline{G}(x(1-\rho)))$ which is negligible with respect to
$\overline{G_I}(x(1-\rho))$. Therefore, one can choose an integer-valued 
$h(x)\to \infty$ such that $\dd{P}\left(A_-^{(0,h(x))}(x)\right)=o(\overline{G_I}(x(1-\rho)))$. So we may apply again the SLLN,  
$\widetilde{t}_{-k} \sim ak$ for sufficiently large $k$, to get
\begin{align*}
 & \bigcup_{k=1}^{\infty} \left(A^{H}_{-k} \bigcap \left(\bigcup_{i=1}^{k} \left\{\tau^{H}_{-k} \ge k+1, \sigma^{H}_{-k+i} > (-\widetilde{t}_{-k}+x)(1-\rho)\right\}\right)\right) \\
 & \quad \simeq \bigcup_{k=1}^{\infty} \bigcup_{\ell=0}^{k-1} \left(A^{H}_{-k} \bigcap \left\{\tau^{H}_{-k} \ge k+1, \sigma^{H}_{-\ell} > (ak+x)(1-\rho)\right\} \right)\\
 & \quad = \bigcup_{\ell=0}^{\infty} \bigcup_{k=\ell+1}^{\infty} \left(A^{H}_{-k} \bigcap \left\{\tau^{H}_{-k} \ge k+1, \sigma^{H}_{-\ell} > (ak+x)(1-\rho)\right\} \right)\\
& \quad \subseteq 
\bigcup_{\ell=0}^{\infty} \left\{\sigma^{H}_{-\ell} > (a(\ell+1) +x)(1-\rho)\right\}.
\end{align*}

On the other hand, for $h(x)\uparrow\infty$ sufficiently slowly and for an appropriate sequence $\varepsilon_{\ell}\downarrow 0$ (that comes from the SLLN), we have  
\begin{align*}
& \ \  \quad \bigcup_{\ell=0}^{\infty} \left\{\sigma^{H}_{-\ell} > (a(\ell+1) +x)(1-\rho)\right\}
 \sim \bigcup_{\ell=h(x)}^{\infty} \left\{\sigma^{H}_{-\ell} > (a(\ell+1) +x)(1-\rho)\right\}\\
& \sim \bigcup_{\ell=h(x)}^{\infty} \left\{\sigma^{H}_{-\ell} > ((a+2\varepsilon_{\ell})(\ell+1) +x+h(x))(1-\rho)+h(x)\right\}\bigcap D_{-(\ell+1)} \equiv E(x),
\end{align*}
where 
\begin{align*}
D_{-\ell} = \bigcap_{j=1}^{\infty}\left\{\sum_{i=-\ell+1}^{-\ell+j} t_i \le (a+\varepsilon_j)j+h(x), \sum_{i=-\ell+1}^{-\ell+j} \sigma^H_i\ge (b^h-\varepsilon_j)j-h(x)\right\}.
\end{align*}
Since $B^0>x$ on the event $E(x)$, we arrive at the following PSBJ for the stationary busy period.
\begin{lemma}
\label{lem:psbj-st1}
If the $GI/GI/1$ feedback queue is stable and its service time distribution has an $\sr{IRV}$ distribution with a finite mean, then
\begin{align}
\label{eq:psbj-st4}
  \{B^{0} > x\} \simeq \bigcup_{k=0}^{\infty} \left\{\sigma^{H}_{-k} > (x+a(k+1))(1-\rho)\right\}.
\end{align}
\end{lemma}

The lemma implies that 
\begin{align}\label{eq:psbj-st3}
\dd{P}(B^0>x) \sim  
\sum_{k=0}^{\infty} \dd{P}(\sigma_{-k}^H > (x+a(k+1))(1-\rho ))
\end{align}
since the sum of the probabilities of pairwise intersections is of order $O(\overline{G_{I}}^2(x))=o(\overline{G_{I}}(x))$.

Then we may conclude that the principle of a single big jump can be applied to the stationary sojourn time too:
\begin{align}\label{eq:negl-h}
\dd{P} (U^0>x) & \sim \sum_{n=0}^{\infty} \dd{P}(U^0>x, \sigma_{-n}^H >(x+(n+1)a)(1-\rho)) \nonumber\\
& \sim \sum_{n=h(x)}^{\infty} \dd{P}(U^0>x, \sigma_{-n}^H >(x+(n+1)a)(1-\rho))
\end{align}
where the second equivalence is valid for any integer-valued function $h(x)\uparrow \infty$, $h(x)=o(x)$ and follows from \eq{decomp} and from the properties of ${\cal IRV}$ and integrated tail distributions, see \app{properties}.

\subsection{Proof of \thr{stationary case}}
\label{sect:tail-st}

First, we comment that it is easy to obtain the logarithmic asymptotics for the stationary sojourn time.  
Since the sojourn time of the customer entering the stationary queue at time $0$ is not bigger than the stationary busy period
and is not smaller than the stationary sojourn time in the auxiliary queue without feedback, and since both bounds have tail distributions that are proportional to the integrated tail distribution of a single service time (see the Appendix for definitions),
we immediately get the logarithmic tail asymptotics:
\begin{align}\label{eq:log-stat}
\log \dd{P}(U^0>x) \sim \log \overline{H_{I}}(x) \sim \log \overline{G_{I}}(x).
\end{align}

Now we provide highlights for obtaining the exact tail asymptotics for the stationary sojourn time distribution and give the final answer. For this, we use the following simplifications, which are made rigorous in the ``WHP'' terminology and due to the $o(x)$-insensitivity of the service-time distribution.
\begin{itemize}
\item [(1)] We observe that the order of services prior to time $0$ is not important for the customer that enters the stationary queue
at time $0$: the joint distribution of
the residual service time and of the queue length at time 0 stays the same for all reasonable service disciplines (that do not allow processor sharing). 
So we may assume that,
up to time $0$, all arriving customers are served in order of their external arrival: the system
serves the ``oldest'' customer a geometric number of times and then turns to the service of the next customer.
\item [(2)] We simplify the model by assuming that all inter-arrival times are deterministic and equal to $a=\lambda^{-1}$.
\item [(3)] We further assume that all service times of all customers but one are equal to $b$, so every customer
but one has a geometric number of services of length $b$. The ``exceptional'' customer may be any customer
$-n\le 0$, it has a geometric number of services, one of those is random and large and all others equal to $b$.
So the total service time of the ``exceptional'' customer has the tail distribution equivalent to
$$
\dd{E}K \cdot \overline{G}(x) = \frac{1}{q}\overline{G}(x).
$$
\item [(4)] We assume that the ``exceptional'' customer arrives at an empty queue, that is, the workload found by this customer is negligible compared with his exceptional service time.
\end{itemize}

Due to the arguments explained above, we can show that the tail asymptotics of the sojourn time of customer $1$ in the original and in the auxiliary
system are equivalent. We start by repeating our calculations from the proof of \thr{U asym 1}, but in two slightly different settings. 

Assume all service times but the very first one are equal to $b$ for the exceptional customer arriving at or before time $0$. Assume that, if customer $1$ arriving at time $0$ is not exceptional, then it finds $X_0=N$ customers in the queue, and otherwise it finds a negligible number of customers compared  with $N$ while its first service time is $Nb$. Assume customer $1$ leaves the system after $K=k$ services. Denote, as before, by $U_i$ the time between its $(i-1)$st and $i$th services and by $X_i$ the queue behind customer $1$ after its $i$th service completion.
How large should $N$ be for the sojourn
time of customer $1$ to be bigger than $x$ where $x$ is large? 

(A) Assume that the (residual) service time, $z$, of the very first customer in the queue
is not bigger than $b$ (so we may neglect it).
When $N$ is large, we get that $U_{1} \sim Nb$. Then we have
\begin{align*}
 X_{i} \sim X_{i-1} p + \lambda U_{i}, \qquad U_{i} \sim X_{i-1} b, \qquad i =1,2,\ldots,k.
\end{align*}
Hence, $X_i \sim X_{i-1}p + \lambda U_{i} \sim N r^i$ and $U_{i} \sim Nbr^{i-1}$ where $r=p+\lambda b<1$. Then $U_1+\ldots+U_{k}\sim Nb(1-r^k)/(1-r)$. Thus, we may conclude that 
\begin{align}\label{eq:XU}
\dd{P}(U>x, K=k) \sim qp^{k-1}\dd{P}(Nb> x_k), 
\end{align}
where $x_k = x(1-r)/(1-r^k)$.

(B) Assume now that both $X_0=N$ and $z$ are large. Then $U_{1} \sim z+Nb$ and $X_1 \sim Np+\lambda (z+Nb)$ and, further, $X_{i}\sim X_{i-1}r \sim X_1 r^{i-1}$ and
$U_{i+1} \sim X_ib \sim X_1b r^{i-1} \sim Nbr^i + zr^i$. Thus,
 $U_1+\ldots+U_{k}\sim Nb(1-r^k)/(1-r)+ z(1-r^k)/(1-r)$ and
 \begin{align}\label{eq:XU2}
 \dd{P}(U>x, K=k) \sim qp^{k-1}\dd{P}(Nb+z> x(1-r)/(1-r^k)).
\end{align} 

Let $W(t)$ be the total work in the system at time $t$. We illustrate $W(t)$ below to see how the cases (A) and (B) occur.

\begin{figure}[h] 
   \centering
   \includegraphics[height=4.5cm]{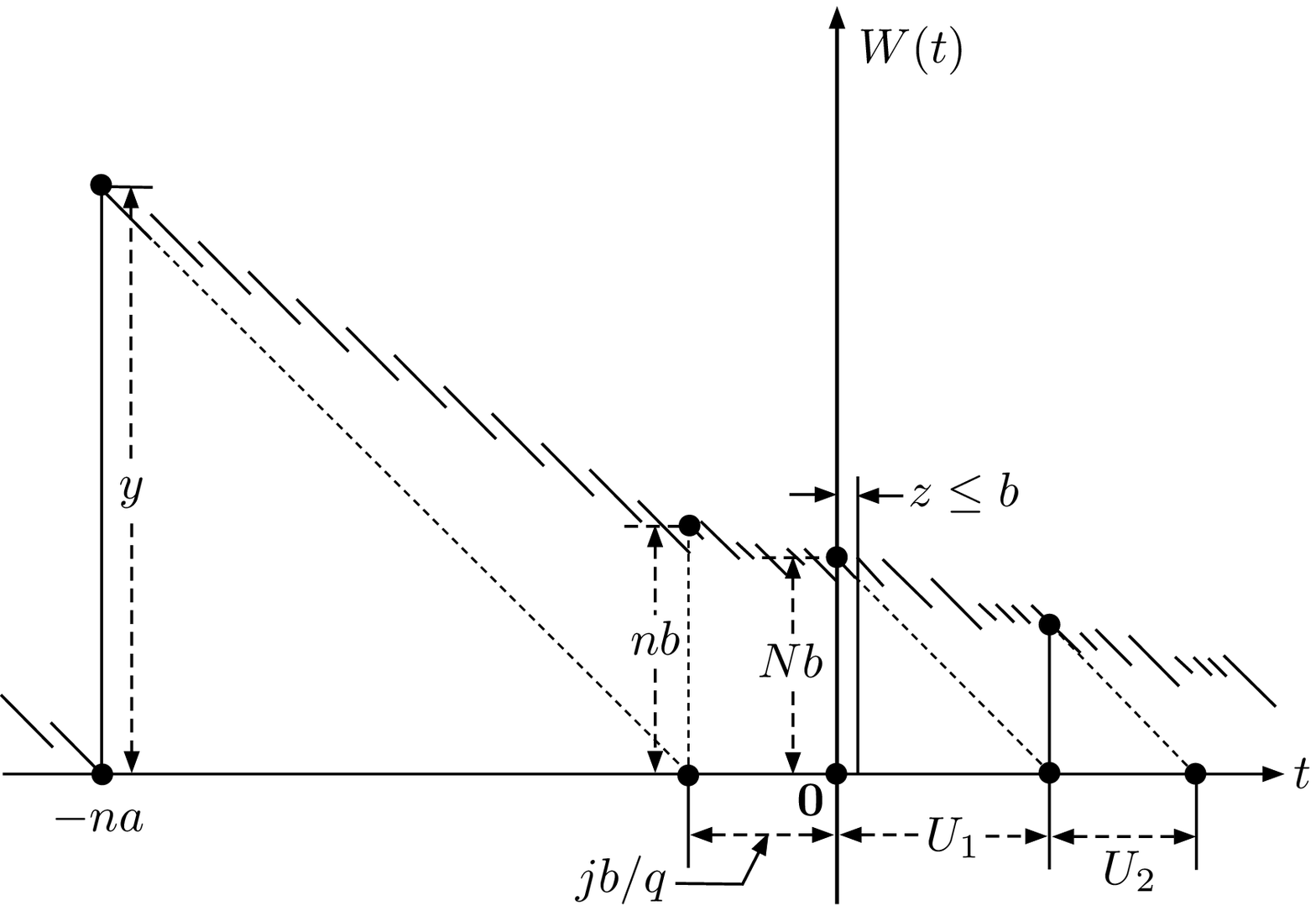} \hspace{1ex}
   \includegraphics[height=4.5cm]{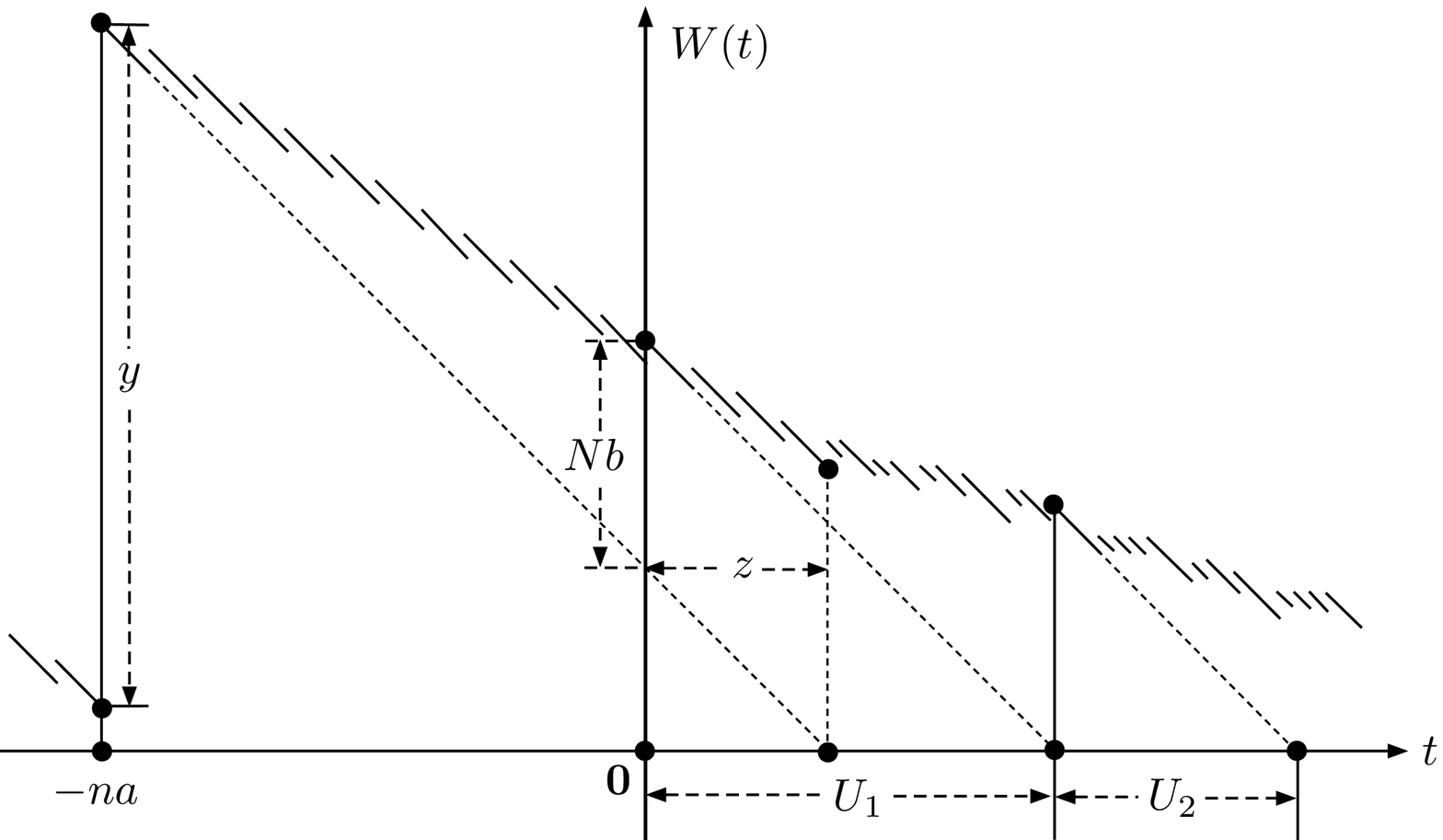}\\
\caption{Sample path of the workload $W(t)$ for the cases (A) and (B) with $K=2$}
\label{fig:Feedback_stationary}
\end{figure}
We will see now that if $K=k$ and if there is a big service time of the $(-n)$th ``exceptional'' 
customer, then the case (A) occurs if $n>x_k/b$ and the case (B) if $n<x_k/b$.

Let the big service time take value $y \gg 1$.  Recall from \eq{negl-h} that it is enough to consider values of $n\ge h(x)$ only, 
where $h(x)\uparrow\infty$, $h(x)=o(x)$. 

For any $k\ge 1$, assume $K=k$ and $y \le na$, then the exceptional service is completed before or at time $0$, and the situation (A) occurs. Hence, $X_{0} \equiv N = n - j$ for some nonnegative $j \le n$, and $y+jb/q \approx na$ because approximately $j$ further customers leave the system prior to time $0$. Then $U \sim Nb (1-r^{k})(1-r)$, and $U > x$ is asymptotically equivalent to $Nb (1-r^{k})(1-r) > x$, where the last inequality is identical with $(n-j)b > x_{k}$. This together with $j \approx (na-y) q/b$ implies that
\begin{align*}
  y \gtrsim x_{k}/q + n(a-b/q).
\end{align*}
Since $na \ge y$, this further implies that $n \gtrsim x_{k}/b$.

We next assume $K=k$ and $n<x_k/b$. Then, the contraposition of the above implication implies that $y > na$, and the situation (B) occurs. Therefore, we should take $y=z+na$, and $U_{0} = nb + z$, where $N=n$. Since $U \sim (nb+z) (1-r^{k})/(1-r)$, $U > x$ is equivalent to $nb + z \gtrsim x_{k}$, and therefore
\begin{align*}
  y \gtrsim x_{k} + n(a-b).
\end{align*}

Combining together both cases, we obtain the following result:
\begin{align}\label{eq:final1-stat}
\dd{P}(U^0>x) & \sim q\sum_{k=1}^{\infty} p^{k-1} 
\Bigg( \sum_{n=1}^{x_k/b - 1} \dd{P}(\sigma>x_k+n(a-b)) \nonumber\\
& \hspace{16ex} +
\sum_{n=x_k/b}^{\infty} \dd{P}(\sigma > x_k/q + n(a-b/q))\Bigg).
\end{align}
Clearly, the second sum in the parentheses is equivalent to
\begin{align*}
\sum_{n=0}^{\infty} \dd{P}(\sigma > x_ka/b + n(a-b/q)) & \sim
(a-b/q)^{-1} \overline{G_{I}}(x_ka/b)
\end{align*}
while the first sum in the parentheses is 
\begin{align*}
\sum_{n=1}^{\infty} - \sum_{n=x_k/b}^{\infty} & \sim
(a-b)^{-1} \left(\overline{G_{I}}(x_k) - \overline{G_{I}}(x_ka/b)\right).
\end{align*}
Hence, we have
\begin{align}
\label{eq:final2-stat}
\dd{P}(U^0>x) & \sim \frac{q}{a-b}\sum_{k=1}^{\infty} p^{k-1} 
\left(\overline{G_{I}}(x_k) + \frac{b(1-q)}{aq-b}\overline{G_{I}}(x_ka/b)
\right),
\end{align}
where we recall that $x_k=\frac{x(1-r)} {1-r^k}$, for $k\ge 1$. Since $m^{(0)}_{K} = \dd{E}\left( \frac {1-r^{K}} {1-r}\right)$, 
we arrive at \eq{final2-stat 2} and \eq{final3-stat}.

\appendix

\section*{Appendix}
\label{sect:appendix}
\setnewcounter 

\section{Properties of heavy-tailed distributions}
\label{app:properties}

We revise basic properties of several classes of heavy-tailed distributions listed at the end of \sectn{introduction} (see \cite{FKZ2013} for the modern theory of heavy tailed distribution and other books \cite{As2003}, \cite{EKM1997} for more detail), and formulate a part of the main result from \cite{FoZa2003} that plays
an important role in our analysis. 

Let $\{ \xi_n\}_{-\infty}^{\infty}$ be $i.i.d.$ r.v.'s
with finite mean $\dd{E} \xi_1$ and with $\dd{P} (\xi_1>0)>0$ and $\dd{P} (\xi_1<0)>0$. 
Let $F(x) = \dd{P} (\xi_1 \leq x)$ be their common distribution and $\overline{F}(x) =
1-F(x) = \dd{P} (\xi_1>x)$ its tail. 
Let $m^+ \equiv m^+(F) = \dd{E} \max (0,\xi_1) = \int_0^{\infty} \overline{F}(t) dt$. 
Let $S_0 =0$, $S_n = \sum_1^n \xi_i$, and 
$\tau = \min \{ n \ge 1 ~:~ S_n \leq 0\}$. 

If $F \in \sr{L}$ (long tailed), that is, \eq{long} holds for some $y>0$, then it holds for all $y$ and, moreover, uniformly in $|y|\le C$, for any fixed $C$.
Therefore, if $F\in {\cal L}$, then there exists a positive function $h(x)\to\infty$ such
that $\overline{F}(x-h(x))\sim \overline{F}(x) \sim \overline{F}(x+h(x))$. In this case we say that
the tail distribution $\overline{F}$ {\it is} $h$-{\it insensitive}.

In what follows, we make use of the following characteristic result (see Theorem 2.47 in \cite{FKZ2013}):
\begin{equation}
\label{charIRV}
F\in {\cal IRV} \ \mbox{if and only if} \ F \ \mbox{is} \ h-\mbox{insensitive, for any} \ h(x)=o(x).
\end{equation}
In particular, if $F$ is $h$-insensitive, then $F$ is $h_c$-insensitive for any $c>0$, where $h_c(x)=h(cx).$

We also use another characteristic result which is a straightforward minor   extension of Theorem 2.48 from \cite{FKZ2013}:
\begin{equation}
\label{char2IRV}
F\in {\cal IRV} \  \ \mbox{if and only if} \ \ \overline{F}(V_n)\sim \overline{F}(v_n),
\end{equation}
for any sequence of non-negative random variables $V_n$ with corresponding means $v_n=\dd{E} V_n$ satisfying
\begin{equation*}
V_n\to\infty \ \mbox{and} \ V_n/v_n\to 1 \ \ \ \mbox{in probability}. 
\end{equation*}

Here is another good property of ${\cal IRV}$ distributions. Let random variables $X$ and $Y$ have arbitrary joint distribution, with the distribution of $X$ being ${\cal IRV}$ and 
$\dd{P} (|Y|>x) =o(\dd{P}(X>x))$. Then
\begin{equation}\label{o-little}
\dd{P} (X+Y>x) \sim \dd{P} (X>x) \ \mbox{as} \ x\to\infty.
\end{equation}

If $F$ is an $\sr{IRV}$ distribution with finite mean, then the distribution with the integrated tail
$\overline{F_I}(x)=\min (1,\int_x^{\infty}\overline{F}(y) dy)$ is also $\sr{IRV}$ and $\overline{F}(x) = o(\overline{F_I}(x))$ and, moreover, $\int_x^{x+h(x)}\overline{F}(y) dy
= o(\overline{F_I}(x))$ if $F_I$ is $h$-insensitive.

We use the following well-known  result: if $\{\sigma_{1,j}\}$ is an $i.i.d.$  sequence of random
variables with common subexponential distribution $F$ and if the counting random variable $K$ does
not depend on the sequence and has a light-tailed distribution, then
\begin{equation}\label{Kesten}
\left\{ \sum_1^K \sigma_{1,j} >x \right\} \simeq \cup_{1}^K\{\sigma_{1,j}>x\} \ \ \mbox{and} \ \ 
\dd{P} \left( \sum_1^K \sigma_{1,j} >x\right) \sim \dd{E} K \overline{F}(x), \quad x\to\infty.
\end{equation}
Here is the principle of a single big jump again: the sum is large when one of the summands is large.

Let $M =\sup_{n\ge 0} \sum_{i=1}^n \xi_i$ where $\{\xi_i\}$ are $i.i.d.$  r.v.'s with negative mean $-m$ 
and with common distribution function $F$ such that $F_I$ is subexponential. Then
\begin{align}
\label{MRSBJ}
 & \{M>x\} \simeq \cup_{n\ge 1}\{\xi_n> x+mn\}, \nonumber\\
 & \dd{P}(M>x) \sim \sum_{n\ge 1} \dd{P}(\xi_1>x+mn) \sim \frac{1}{m}\overline{F_I}(x).
\end{align} 
Further, if $F_I$ is subexponential, then, for any sequence $m_n\to m>0$ and any function
$h(x)=o(x)$,
\begin{align}\label{PSBJ100}
\sum_{n\ge 0} \overline{F}(x+m_nn+h(x))\sim \frac{1}{m}\overline{F_I}(x)
\end{align}
and, for any sequence $c_n\to 0$,
\begin{align}\label{PSBJ101}
\sum_{n\ge 0} c_n\overline{F}(x+m_nn+h(x)) =o\left(\overline{F_I}(x)\right).
\end{align}

\begin{remark}
\label{rem:tailprop}
Let ${\cal K}$ be any of the classes ${\cal L},{\cal RV}, {\cal IRV},
{\cal  D},
 {\cal S}, {\cal S}^*$.
The property of belonging to class ${\cal K}$ is
a {\it tail property}: if $F\in {\cal K}$ and if $\overline{G}(x) \sim C\overline{F}(x)$ where $C$ is a positive constant, then
$G\in {\cal K}$. In particular,\\
(i) if $F\in {\cal K}$, then $F_+\in {\cal K}$;\\
(ii) if the random variable $\xi$ has distribution $F\in {\cal K}$ and $c_1 > 0$ and $c_2$ are any constants,
then the distribution of the random variable $\eta = c_1\xi + c_2$ also belongs to ${\cal K}$; \\
(iii) if the random variable $\xi$ may be represented as $\xi = \sigma - t$ where $\sigma$ and $t$ are mutually independent random variables and $t$ is non-negative (or, slightly more generally, bounded from below), and if the distribution of
$\sigma$ belongs to class ${\cal K}$, then $\dd{P} (\xi >x) \sim
\dd{P} (\sigma > x)$, so 
the distribution of $\xi$ belongs to ${\cal K}$ too.
\end{remark}

The following result is a part of Theorem 1 in \cite{FoZa2003}, see also
\cite{FoPaZa2005} for a more general statement.

\begin{theorem}
\label{thr:FZ1} 
Let $S_n=\sum_1^n \xi$, $S_0=0$ be a random walk with $i.i.d.$  increments with distribution function $F$ and finite negative mean
\begin{equation}\label{NEG}
 \dd{E} \xi_1 = -m < 0.
\end{equation} 
Assume $F\in {\cal S}^*$. Let $T \le \infty$ be any stopping time (with respect to
$\{ \xi_n\}$). Let $M_{T} = \max_{0\le n\le T} S_n$. Then
\begin{equation}\label{si}
 \lim_{x\to\infty} \frac{\dd{P} (M_{T} >x)}{\overline{F}(x)} = \dd{E} T.
\end{equation}
\end{theorem}

\section {Proof of \cor{gg2}}
\label{app:gg2}

We first note that the event $\{\tau^{H} \ge n \}$ is independent of $\sigma^{H}_{n}$ and $\xi^{H}_{n}$ because $\{\tau^{H} \ge n \} = \{\tau^{H} \le n-1 \}^{c}$ is $\sigma(\{\xi^{H}_{\ell}; 1 \le \ell \le n-1\})$-measurable. Hence, 
\begin{align*}
  \dd{P}(\tau^{H} \ge n, \sigma^{H}_{n} > x(1-\rho)) = \dd{P}(\tau^{H} \ge n) \dd{P}(\sigma^{H}_{n} > x(1-\rho)),
\end{align*}
and therefore the equivalence \eq{PSBJ-B1} is immediate from \eq{B3} of \thr{gg1}, while \eq{PSBJ-B2} easily follows from \eq{PSBJ-B1}.

Thus, it remains to prove \eq{PSBJ-B3}. For this, we introduce some notation. Let $S^{H}_{n} = \sum_{i=1}^n \xi^{H}_{i}$ and let $S^{\sigma^{H}}_n = \sum_{i=1}^n \sigma^{H}_i$. Define a sequence of events $E_{n}$, $n=0,1,\ldots$, as
 $$
 E_{n} = \cap_{\ell\ge 1}\{(S^{H}_{\ell+n}-S^{H}_{n})\ge (b_{H}-a){\ell}-\delta_{\ell} \ell - C, 
 (S^{\sigma^{H}}_{\ell+n}-S^{\sigma^{H}}_{n})\ge b_{H}{\ell} - \delta_{\ell} \ell -C\}
 $$
 which is stationary in $n$ (here, by convention, $S^{H}_0=S^{\sigma^{H}}_0=0$). Due to the SLLN,
there exists a sequence $\delta_{\ell} \downarrow 0$ such that
$$
\dd{P} (|S^{H}_{\ell}/\ell-(b_{H}-a)|\le \delta_{\ell} \ \mbox{and} \ 
 |S^{\sigma^{H}}_{\ell}/{\ell} - b_{H}|\le \delta_{\ell}, \ \mbox{for all} \ {\ell}\ge n)\to 1,
 $$
 as $n\to\infty$. Therefore, for any $\varepsilon >0$, there exists $C=C_{\varepsilon}>0$, for which $E_{n}$ is denoted by $E_{n,\varepsilon}$, such that
\begin{align}
\label{eq:En 1}
 \dd{P}(E_{n,\varepsilon}) = \dd{P}(E_{0,\varepsilon}) \ge 1- \varepsilon.
\end{align}
 Introduce a function $h_{\varepsilon}(x)$ by
 $$
 h_{\varepsilon}(x)= \max_{i \le [x/a]} \left(i \delta_{i}\right) + C_{\varepsilon} + b_{H},
 $$
 where $[x/a]$ is the integer part of the ratio $x/a$. Then, for this $\varepsilon$ and $n \ge 1$, define $J_{n,\varepsilon}(x)$ as
\begin{align*}
  J_{n,\varepsilon}(x) = \left\{ \tau^{H} \ge n, \xi^{H}_{n} > x(1-\rho) + h_{\varepsilon}(x)\right\}, \qquad x > 0.
\end{align*}
 
 Then, on the event $J_{n,\varepsilon}(x) \cap E_{n,\varepsilon}$, we have $S^{H}_{n-1} > 0$, $S^{H}_{n} > \xi^{H}_{n} > x(1-\rho) + h_{\varepsilon}(x)$, and therefore
\begin{align*}
  S^{H}_{n+\ell} & > x(1-\rho) + h_{\varepsilon}(x) - (1-\rho) a \ell - \delta_{\ell} \ell - C_{\varepsilon}\\
  & = (1-\rho) (x-a \ell) + \max_{i \le [x/a]} \left(i \delta_{i}\right) - \delta_{\ell} \ell + b_{H} > 0, \qquad 0 \le \ell \le [x/a].
\end{align*}
Hence, letting $\ell_{0} = [x/a]$, we have, on the same event,
\begin{align*}
  B \ge \sum_{i = n}^{n + \ell_{0}} \sigma^{H}_{i} & = \sigma^{H}_{n} + S^{\sigma^{H}}_{n+\ell_{0}} - S^{\sigma^{H}}_{n} > x(1-\rho) + h_{\varepsilon}(x) + b_{H} \ell_{0} - \delta_{\ell_{0}} \ell_{0} - C_{\varepsilon} > x.
\end{align*}
For any integer $N \ge 1$, let
\begin{align*}
  L_{N,\varepsilon}(x) = \cup_{n=1}^N \{\max_{i<n} \xi^{H}_i \le x(1-\rho)\} \cap J_{n,\varepsilon}(x) \cap E_{n,\varepsilon},
\end{align*}
then we have
\begin{align}
\label{eq:DN 1}
  \{B>x\} \supseteq L_{N,\varepsilon}(x).
\end{align}

 Let $F^{H}$ be the distribution of $\xi^{H}$, and recall that $H$ is the distribution of $\sigma^{H}$. Both of them are intermediate regularly varying, they are tail-equivalent and $h$-insensitive (see \rem{tailprop} and \eqref{charIRV}). Since $L_{N,\varepsilon}(x)$ is a union of $N$ disjoint events and $\{\tau^{H} \ge n, \max_{i<n} \xi^{H}_i \le x(1-\rho) \}$ is independent of $\xi^{H}_{n}$ and $E_{n,\varepsilon}$, \eq{En 1} yields, as $x\to\infty$, 
\begin{align}
\label{eq:LN 1}
   \dd{P} (L_{N,\varepsilon}(x)) 
 &= \sum_{n=1}^N \dd{P} (\tau^{H} \ge n, \max_{i<n} \xi^{H}_i \le x(1-\rho)) \cdot 
 \overline{F^{H}}( x(1-\rho ) + h_{\varepsilon}(x))\cdot \dd{P} (E_{n,\varepsilon}) \nonumber\\
 &\sim
 \sum_{n=1}^N \dd{P} (\tau^{H} \ge n) \overline{H}(x(1-\rho )) \dd{P} (E_{n,\varepsilon}) \nonumber\\
 &\ge
 (1-\varepsilon ) \overline{H}(x(1-\rho )) \sum_{n=1}^N \dd{P} (\tau^{H} \ge n).
\end{align}
Let $L_{\varepsilon}(x) = \lim_{N \to \infty} L_{N,\varepsilon}(x)$, then 
$L_{\varepsilon}(x) \subset \{B>x\}$ by \eq{DN 1}, and, for any $N$, by \eq{B3} of \thr{gg1},
$$
\dd{P} (B>x) - \dd{P} (L_{\varepsilon}(x)) \lesssim  
 {\overline H}(x(1-\rho))\left(\varepsilon \sum_{n=1}^N \dd{P}(\tau^{H} \ge n) + \sum_{n=N+1}^{\infty}
\dd{P} (\tau^{H} >n)\right).
$$
Choosing $N$ such that $\sum_{n=N+1}^{\infty} \dd{P} (\tau^{H} >n) \le \varepsilon \dd{E} \tau^{H}$, we get
\begin{align}
\label{eq:BL 1}
 0 \le \dd{P} (B>x) - \dd{P} (L_{\varepsilon}(x)) \lesssim  
2 {\overline H}(x(1-\rho))\varepsilon \dd{E} \tau^{H}.
\end{align}

For $x > 0$, define events $J(x)$ and $\ol{J}_{\varepsilon}(x)$ as
\begin{align*}
 & J(x) = \bigcup_{n=1}^{\infty} \{\tau^{H} \ge n, \xi^{H}_{n} > x(1-\rho)\},\\
 & \ol{J}_{\varepsilon}(x) = \bigcup_{n=1}^{\infty} \{\max_{i<n} \xi^{H}_i \le x(1-\rho)\} \cap J_{n,\varepsilon}(x).
\end{align*}
Since $L_{\varepsilon}(x) = \cup_{n=1}^{\infty} \{\max_{i<n} \xi^{H}_i \le x(1-\rho)\} \cap J_{n,\varepsilon}(x) \cap E_{n,\varepsilon}$, we have
\begin{align*}
 0 \le  \dd{P} (\ol{J}_{\varepsilon}(x)) - \dd{P} (L_{\varepsilon}(x)) & = \sum_{n=1}^{\infty} \dd{P}(\tau^{H} \ge n, \max_{i<n} \xi^{H}_i \le x(1-\rho)) \overline{F^{H}}(x(1-\rho)) (1- \dd{P}(E_{n,\varepsilon}))\\
 & \le
\varepsilon \dd{E}\tau^{H} \overline{F^{H}}(x(1-\rho)).
\end{align*}
Further, $\ol{J}_{\varepsilon}(x) \subset J(x)$ and
$$
 0 \le \dd{P} (J(x))-\dd{P}(\ol{J}_{\varepsilon}(x)) \le
\dd{E}\tau^{H} \dd{P} (\xi^{H}_1\in (x(1-\rho), x(1-\rho)+h_{\varepsilon}(x)]) =
o(\overline{F^{H}}(x(1-\rho)).
$$
Combining those with \eq{BL 1}, we obtain that
\begin{align*}
  \dd{P} \left(\{B>x\}\setminus J(x)\right) & \le \dd{P}(B>x) - \dd{P}(L_{\varepsilon}(x) \cap J(x))\\
& = \dd{P}(B>x) - \dd{P}(L_{\varepsilon}(x)) \le 2\varepsilon \dd{E}\tau^{H} \overline{F^{H}}(x(1-\rho)),\\
 \dd{P} \left(J(x)\setminus \{B>x\}\right) & \le \dd{P} \left(J(x)\right) - \dd{P}\left(L_{\varepsilon}(x)\right)\\& \lesssim \varepsilon \dd{E}\tau^{H} \overline{F^{H}}(x(1-\rho)) + o(\overline{F^{H}}(x(1-\rho)).
\end{align*}
Since $\varepsilon >0$ may be taken arbitrarily small, we arrive at \eq{PSBJ-B3}.

\section{Alternative proof of \cor{U finite 1}}
\label{app:alternative} 

In this section, we give an alternative proof of \cor{U finite 1}, which is based on the result from \cite{FoZa2003}, not using PSBJ. Instead of it, our basic tools are \thr{FZ1} and the law of large numbers. We also slightly generalise \cor{U finite 1}.

\begin{theorem}
\label{thr:U finite a1}
For the stable $GI/GI/1$ feedback queue, assume that its service time distribution is intermediate regularly varying and has a finite mean. If the first customer arriving at the empty system has an exceptional first service time $\eta$ instead of $\sigma_{1,0}$, that is, $U_{1} = \eta$, such that
\begin{itemize}
\item [(I)] $\eta$ has an intermediate regularly varying distribution with $\dd{E}(\eta) < \infty$;
\item [(II)] $\ds \frac {\dd{P}(\eta > x)} {\dd{P}(\sigma > x)}$ converges to a constant from $[0,\infty]$, as $x \to \infty$;
\item [(III)] $\dd{E}(X^{(0)}_{k-1}) < \infty$ for $k \ge 1$;
\end{itemize}
then, for each $k \ge 1$, as $x \to \infty$,
\begin{align}
\label{eq:U finite a1}
  \dd{P}\left(T_{k} > x\right) \sim \sum_{\ell=0}^{k-2} & \dd{E}\left( 1 + X^{(0)}_{k-\ell-1} \right)  \dd{P}\left((1+ m^{(0)}_{\ell}) \sigma > x \right) + \dd{P}\left((1+ m^{(0)}_{k-1}) \eta > x \right).
\end{align}
\end{theorem}
\begin{remark}
\label{rem:U finite a1}
If $\eta = \sigma_{1,0}$, then the conditions (I)--(III) are satisfied, and this theorem is just \cor{U finite 1}.
\end{remark}

\begin{proof}
Recall the definitions of $X_{\ell-1}, U_{\ell}$ under $\sigma_{1,0} = \eta$. We have $X_{0} = u_{0}= 0$, $U_{1} = \eta$, and, for $\ell \ge 2$,
\begin{align}
\label{eq:X0 1}
 & X^{(0)}_{\ell-1} = N^{e}(T_{\ell-2}+U_{\ell-1}) - N^{e}(T_{\ell-2}) + N^{B}_{\ell-2}(X^{(0)}_{\ell-2}),\\
\label{eq:X0 sum}
 &  \sum_{j=1}^{\ell-1}X^{(0)}_{j} = N^{e}(T_{\ell-1}) + \sum_{j=2}^{\ell-1} N^{B}_{j-1}(X^{(0)}_{j-1}),\\
\label{eq:T0 1}
 & T_{\ell} = \sum_{j=1}^{\ell} U_{j} = \sum_{j}^{\ell} \sum_{i=0}^{X^{(0)}_{j-1}} \sigma_{j,i},
\end{align}
where $X_{j-1}$ and $\sigma_{j, i}$ are independent. Since \eq{U finite a1} is an identity for $k=1$, we assume that $k \ge 2$. We partition the event $\{T_{k} > x\}$ into the following $k$ disjoint sets for each $\vc{y} \equiv( y_{1}, y_{2}, \ldots, y_{k-1}) > \vc{0}$. 
\begin{align*}
 I_{k}^{(\ell)}(\vc{y},x) = \begin{cases}
\{T_{k} > x, U_{j} \le y_{j}, 1 \le j \le k-1\}, & \ell=1,\\
 \{T_{k} > x, U_{j} \le y_{j}, 1 \le j \le k - \ell, U_{k-\ell+1} > y_{k-\ell+1}\}, & 2 \le \ell \le k.
\end{cases}
\end{align*}

We prove that
\begin{align}
\label{eq:Ik 1}
 & \dd{P}\left(I_{k}^{(\ell)}(\vc{y},x)\right) \sim \begin{cases}
\dd{E}\left( 1 + X^{(0)}_{k-\ell} \right) \dd{P}\left((1+ m^{(0)}_{\ell-1}) \sigma > x \right), & 1 \le \ell \le k-1,\\
\dd{P}\left((1+ m^{(0)}_{k-1}) \eta > x \right), & \ell=k,
\end{cases}
\end{align}
as $x \to \infty$ then $y_{1}, \ldots, y_{k-1} \to \infty$. By the assumptions (I)--(III), these asymptotics yield \eq{U finite a1}, and therefore the theorem is obtained.

We prove \eq{Ik 1} deriving upper and lower bounds. We first consider the case that $\ell = 1$. Since $U_{j} \le y_{j}$ for $1 \le j \le k-1$, we have that $T_{j} \le \sum_{j'=1}^{j} y_{j'}$ for $1 \le j \le k-1$, and 
\begin{align*}
  \sum_{j=1}^{\ell-1}X^{(0)}_{j} \le N^{e}\left(\sum_{j'=1}^{\ell-1} y_{j'}\right) + \sum_{j=2}^{\ell-1} N^{B}_{j-1}\left(X^{(0)}_{j-1}\right) \quad \mbox{ on } I_{k}^{(1)}(\vc{y},x).
\end{align*}
This inductively shows that $X^{(0)}_{k-1}$ has light tail on $I_{k}^{(\ell)}(\vc{y},x)$, and therefore \eq{T0 1} implies that
\begin{align*}
  \limsup_{x \to \infty} \frac {\dd{P}\left( I_{k}^{(1)}(\vc{y},x) \right)} {\dd{P}(\sigma > x)} \le \limsup_{x \to \infty} \frac {\dd{P}\left( \sum_{i=0}^{X^{(0)}_{k-1}} \sigma_{k,i} > x - \sum_{j=1}^{k-1} y_{j}\right)} {\dd{P}(\sigma > x)} = \dd{E}\left( 1 + X^{(0)}_{k-1} \right).
\end{align*}
The corresponding lower bound is obvious. That is,
\begin{align*}
  \liminf_{x \to \infty} \frac {\dd{P}\left( I_{k}^{(1)}(\vc{y},x) \right)} {\dd{P}(\sigma > x)} & \ge \liminf_{x \to \infty} \frac {\dd{P}\left( \sum_{i=0}^{X^{(0)}_{k-1}1(U_{j} <y_{j}, 1\le j \le k-1)} \sigma_{k,i} > x\right)} {\dd{P}(\sigma > x)}\\
  & = \dd{E}\left( 1(U_{j} <y_{j}, 1\le j \le k-1) \left( 1 + X^{(0)}_{k-1} \right) \right).
\end{align*}
Hence, letting $y_{j} \to \infty$ for $j=1,2,\ldots,k$, we obtain \eq{Ik 1} for $\ell=1$.

We next consider the case $\ell=2$. Let $c_{1}$ be a positive constant, which will be appropriately determined for each sufficiently small $\varepsilon > 0$.
\begin{align*}
  \dd{P}\left( I_{k}^{(2)}(\vc{y},x) \right) \le A^{(2)}_{+}(y_{1},\ldots,y_{k-2},x) + B^{(2)}_{+}(y_{1},\ldots,y_{k-2},x),
\end{align*}
where
\begin{align*}
 A^{(2)}_{+}(y_{1},\ldots,y_{k-2},x) & = \dd{P}\Bigg( U_{j} \le y_{j}, 1 \le j \le k-2, U_{k-1} > c_{1} x\Bigg),\\
 B^{(2)}_{+}(y_{1},\ldots,y_{k-2},x) & = \dd{P}\Bigg( \sum_{i=0}^{X^{(0)}_{k-1}} \sigma_{k,i} > (1-c_{1}) x - (y_{1}+\ldots+y_{k-2}),\\
 & \hspace{10ex} U_{j} \le y_{j}, 1 \le j \le k-2, y_{k-1} < U_{k-1} \le c_{1} x\Bigg). 
\end{align*}
Since $U_{k-1} = \sum_{i=1}^{X^{(0)}_{k-2}} \sigma_{k-1,i}$, the term $A^{(2)}_{+}(y_{1},\ldots,y_{k-2},x)$ has the same asymptotics as  $I^{(1)}_{k-1}(\vc{y},c_{1}x)$. Hence, it follows from the asymptotics of $I^{(1)}_{k}(\vc{y},x)$ that
\begin{align}
\label{eq:A2 1}
  A^{(2)}_{+}(y_{1},\ldots,y_{k-2},x) \sim \dd{E}\left( 1 + X^{(0)}_{k-1} \right) 1(k=2) \dd{P}(\eta > c_{1}x) + 1(k \ge 3) \dd{P}(\sigma > c_{1}x).
\end{align}
Thus, if we show that $B^{(2)}_{+}(y_{1},\ldots,y_{k-2},x)$ is asymptotically negligible, then the right-hand side of \eq{Ik 1} is obtained as an upper bound for $\ell=2$. To see this, we consider $X^{(0)}_{k-1}$ on the event 
\begin{align*}
  E^{U}_{k-1}(\vc{y},x) \equiv \{U_{j} \le y_{j}, 1 \le j \le k-2, y_{k-1} < U_{k-1} \le c_{1} x\},
\end{align*}
on which we have
\begin{align*}
  X^{(0)}_{k-1} = N^{e}(T_{k-2}+U_{k-1}) - N^{e}(T_{k-2}) + N^{B}_{k-1}(X^{(0)}_{k-2}) \lesssim N^{e}(c_{1}x), \qquad x \to \infty.
\end{align*}
Hence, letting
\begin{align*}
 & z_{k-2} = y_{1}+\ldots+y_{k-2},\\
 & \widetilde{S}^{\sigma}_{n} = \sum_{i=1}^{n} (\sigma_{k,i} - (b+\varepsilon)), \qquad \widetilde{M}_{n} = \max_{1 \le j \le n} \widetilde{S}^{\sigma}_{j}, \qquad n \ge 1,
\end{align*}
and applying \thr{FZ1}, we have
\begin{align*}
  B^{(2)}_{+}(y_{1},\ldots,y_{k-2},x) & = \dd{P}\Big( \widetilde{S}^{\sigma}_{X^{(0)}_{k-1}+1} + (b+\varepsilon) (X^{(0)}_{k-1}+1) > (1-c_{1}) x - z_{k-2}, E^{U}_{k-1}(\vc{y},x) \Big)\\
 & \le \dd{P}\Big( \widetilde{M}_{X^{(0)}_{k-1}+1} + (b+\varepsilon) (X^{(0)}_{k-1}+1) > (1-c_{1}) x - z_{k-2}, E^{U}_{k-1}(\vc{y},x) \Big)\\
 & \le \dd{P}\Big( \widetilde{M}_{X^{(0)}_{k-1}+1} > (1-c_{1}(1+(b+\varepsilon)(\lambda+\varepsilon)) x, E^{U}_{k-1}(\vc{y},x) \Big)\\
 & \quad + \dd{P}\Big( X^{(0)}_{k-1}+1 > c_{1} (\lambda+\varepsilon)x - z_{k-2}, E^{U}_{k-1}(\vc{y},x) \Big)\\
 & \lesssim \dd{E}\left( (X^{(0)}_{k-1}+1)1_{E^{U}_{k-1}(\vc{y},x)} \right) \dd{P}\left( \sigma > (1-c_{1}(1+(b+\varepsilon)(\lambda+\varepsilon)) x \right)\\
 & \quad + \dd{P}\Big( N^{e}(c_{1}x) > c_{1} (\lambda+\varepsilon)x - z_{k-2} \Big),
\end{align*}
where the last probability term decays super-exponentially fast, so it is negligible. Thus, if  we choose $c_{1} > 0$ such that $c_{1}(1+(b+\varepsilon)(\lambda+\varepsilon)) < 1$, then $B^{(2)}_{+}(y_{1},\ldots,y_{k-2},x)$ is asymptotically negligible because
\begin{align*}
  \lim_{x \to \infty} \dd{E}\left( (X^{(0)}_{k-1}+1)1_{E^{U}_{k-1}(\vc{y},x)} \right) = 0.
\end{align*}
Consequently, we choose $c_{1} = \frac {1 - \varepsilon}{1 + (b+\varepsilon)(\lambda+\varepsilon)}$, which converges to $(1 + \lambda b)^{-1}$ as $\varepsilon \downarrow 0$. Thus, we have proved that the right-hand side of \eq{Ik 1} is an upper bound for $\ell=2$.

For the lower bound for $\ell=2$, we take another decomposition. Let $d_{1} = \frac {1+\varepsilon}{1+\lambda b} < 1$ for a sufficiently small $\varepsilon > 0$, then, for $d_{1} x > y_{k-1}$,
\begin{align*}
  \dd{P}\left( I_{k}^{(2)}(y_{1},\ldots,y_{k-2},d_{1}x,x) \right) & \ge \dd{P}(U_{k-1} > d_{1} x, U_{j} \le y_{j}, 1 \le j \le k - 2)\\
  & \qquad - \dd{P}(T_{k} \le x, U_{j} \le y_{j}, 1 \le j \le k - 2, U_{k-1} > d_{1} x).
\end{align*}
Similar to \eq{A2 1}, we have
\begin{align*}
 & \dd{P}(U_{k-1} > d_{1} x, U_{j} \le y_{j}, 1 \le j \le k - 2) \\
 & \quad \sim \dd{E}\left( 1 + X^{(0)}_{k-1} \right) 1(k=2) \dd{P}(\eta > d_{1}x) + 1(k \ge 3) \dd{P}(\sigma > d_{1}x).
\end{align*}
On the other hand, by the law of large numbers,
\begin{align*}
  & \dd{P}(T_{k} \le x, U_{j} \le y_{j}, 1 \le j \le k - 2, U_{k-1} > d_{1}x)\\
  & \quad \le \dd{P}\left(\sum_{i=0}^{X^{(0)}_{k-1}} \sigma_{k,i} \le (1-d_{1} x), U_{j} \le y_{j}, 1 \le j \le k - 2, U_{k-1} > d_{1} x\right)\\
  & \quad \le \dd{P}\left(\sum_{i=0}^{N^{e}(z_{k-2}+d_{1}x)-N^{e}(z_{k-2}) + N^{B}_{k-1}(X^{(0)}_{k-2})} \hspace{-10ex} \sigma_{k,i} \le (1-d_{1} x), U_{j} \le y_{j}, 1 \le j \le k - 2, U_{k-1} > d_{1} x\right)\\
  & \quad = o(1) \dd{P}\left(U_{k-1} > d_{1} x\right).
\end{align*}
Hence, this term is asymptotically negligible, and therefore we have the asymptotic lower bound for $I_{k}^{(2)}(\vc{y},x)$, which agrees with the upper bound, by letting $\varepsilon \downarrow 0$. Thus, we have proved \eq{Ik 1} for $\ell=2$. For $\ell = 3,\ldots,k$, \eq{Ik 1} is similarly proved (we omit the details). Then the proof  of the corollary is completed.
\end{proof}

\subsection*{Acknowledgements}
The authors are grateful to the reviewers for helpful comments and suggestions.


\begin{thebibliography}{99}

\bibitem[{Asmussen(2003)}]{As2003}
\textsc{Asmussen, S.}(2003).
\newblock \textit{Applied Probability and Queues.} 
\newblock Springer, 2nd Edition.

\bibitem[{Baccelli-Foss(2004)}]{BaFo2004}
\textsc{Baccelli, F. and Foss, S.} (2004).
\newblock Moments and tails in monotone-separable stochastic networks.
\newblock \textit{Annals of Applied Probability}, \textbf{14}, 612--650.

\bibitem[{EmKluMik(1997)}]{EKM1997}
\textsc{Embrechts, P.; Kl{\"u}ppelberg, C.; Mikosch, T.} (1997).
\newblock Modelling Extremal Events.
\newblock \textit{Springer}. 

\bibitem[{Foss-Korshunov(2012)}]{FoKo2012}
\textsc{Foss, S.; Korshunov, D.} (2012).
\newblock On Large Delays in Multi-Server Queues with Heavy Tails.
\newblock \textit{Mathematics of Operations Research}, \textbf{37}, 201--218.

\bibitem[{Foss-Korshunov-Zachary(2013)}]{FKZ2013}
\textsc{Foss, S.; Korshunov, D.; Zachary, S.} (2013).
\newblock \textit{An Introduction to Heavy-Tailed and Subexponential Distributions}.
\newblock Springer, 2nd Edition. 

\bibitem[{FosMiy(2014)}]{FoMy2014}
\textsc{Foss, S.; Miyazawa, M.} (2014).
\newblock Two-node fluid network with a heavy-tailed random input: the strong stability case.
\newblock \textit{Journal of Applied Probability}, \textbf{51A}, 249--265.


\bibitem[{FosPalZac(2005)}]{FoPaZa2005}
\textsc{Foss, S.; Palmowski, Z.; Zachary, S.}(2005)
\newblock The probability of exceeding a high boundary on a random time interval for a heavy-tailed random walk.
\newblock \textit{Annals of Applied Probability}, \textbf{15}, 1936--1957.

\bibitem[{FosPuh(2011)}]{FoPu(2011)}
\textsc{Foss, S.; Puhalskii, A.} (2011).
\newblock On the limit law of a random walk conditioned to reach a high level.
\newblock \textit{Stochastic Processes and their Applications}, \textbf{121}, 288--313.

\bibitem[{Foss-Zachary(2003)}]{FoZa2003}
\textsc{Foss, S.; Zachary, S.} (2003).
\newblock The maximum on a random time interval of a random walk with a
long-tailed increments and negative drift.
\newblock \textit{Annals of Applied Probability}, \textbf{13}, 37--53.

\bibitem[{Greenwood(1973)}]{Gree1973}
\textsc{Greenwood, P.} (1973).
\newblock Asymptotics of randomly stopped sequences with independent increments.
\newblock \textit{Annals of Probability}, \textbf{1}, 317--321.

\bibitem[{Jele-Momc-Zwar(2004)}]{JeMoZw2004}
\textsc{Jelenkovic, P., Momcilovic, P. and Zwart, B.} (2004).
\newblock Reduced load equivalence under subexponentiality.
\newblock \textit{Queueing Systems}, \textbf{46}, 97--112. 


\bibitem[{Takacs(1963)}]{Taka1963}
\textsc{Takacs, L.} (1963).
\newblock A single-server queue with feedback.
\newblock \textit{Bell System Technical Journal}, \textbf{42} 505--519.

\bibitem[{Zwart(2001)}]{Zwar2001}
\textsc{Zwart, A.P.} (2001).
\newblock Queueing systems with heavy tails. 
\newblock PhD Thesis. Eindhoven: Technische Universiteit Eindhoven.

\end{thebibliography}

\end{document}